\documentclass{tac}
\pdfoutput=1
\usepackage[toc,page]{appendix}
\usepackage{mathtools} % Loads and extends amsmath
\usepackage{amssymb} % Extra mathematical symbols (loads amsfonts)
\usepackage{stmaryrd} % Some maths symbols for logic and computer science
\usepackage{relsize} % Additional relative sizes for fonts
\usepackage{microtype} % Improves appearance of writing
\usepackage{multicol} % Multi-column environments
\usepackage{csquotes} % Environments for quotes
\usepackage{xspace}
\usepackage{hyperref} % Hyperlink citations, comment out for arxiv submission

\hypersetup{ bookmarksnumbered=true
	   , bookmarksopen=true
	   , breaklinks=true
	   , urlcolor=black
	   , colorlinks%
	   , urlcolor=black
	   , linkcolor=black
	   , citecolor=black
           }
\usepackage[nocompress]{cite} % Comment out if using natbib or apacite
\usepackage[hyperpageref]{backref}
\renewcommand*{\backref}[1]{}
\renewcommand*{\backrefalt}[4]{%
	\ifcase #1 (Not cited.)%
	\or        (Cited on page~#2.)%
	\else      (Cited on pages~#2.)%
  \fi}
\usepackage{etoolbox}
\apptocmd{\sloppy}{\hbadness 10000\relax}{}{}
  \usepackage{graphicx} % Needed to include images
  \usepackage{subcaption} % Needed to define subfigures
  \usepackage[dvipsnames]{xcolor} % Needed to use color names in hyperref options

    \usepackage{tikz} % TikZ
    \usetikzlibrary{
      petri, % To draw petri nets
      backgrounds, % To define image backgrounds
      arrows, % To use and define further arrow tips
      positioning, % To use expressions like "right = 1 of 1"
      decorations.markings, % Needed to define oriented wiring diagrams
      calc,  % Needed to define oriented wiring diagrams
      fit, % Needed to compose wiring diagrams
    }

\usepackage{tabularx, cellspace} % Needed to typeset the tables for listing symbols
\usepackage{tikz-cd}

\usepackage{adjustbox}

  \newcounter{theoremUnified} % Unified coutner for all theorem environments
  \def\thetheoremUnified{\arabic{section}} % Needed to have counters going with sections
  \numberwithin{theoremUnified}{section} % Numbering within sections
  \numberwithin{theoremUnified}{section} % Equations are also numbered within sections

  \theoremstyle{plainStyle} % Environments in plain style
    \newtheorem{example}[theoremUnified]{Example}
      \AfterEndEnvironment{example}{\noindent\ignorespaces}
    
      \AfterEndEnvironment{fact}{\noindent\ignorespaces}
    
      \AfterEndEnvironment{examplehard}{\noindent\ignorespaces}
    \newtheorem{remark}[theoremUnified]{Remark}
      \AfterEndEnvironment{remark}{\noindent\ignorespaces}
    
      \AfterEndEnvironment{remarkhard}{\noindent\ignorespaces}
  \theoremstyle{italicStyle} % Environments in italic style
    \newtheorem{definition}[theoremUnified]{Definition}
      \AfterEndEnvironment{definition}{\noindent\ignorespaces}
    \newtheorem{proposition}[theoremUnified]{Proposition}
      \AfterEndEnvironment{proposition}{\noindent\ignorespaces}
    \newtheorem{lemma}[theoremUnified]{Lemma}
      \AfterEndEnvironment{lemma}{\noindent\ignorespaces}
    
      \AfterEndEnvironment{theorem}{\noindent\ignorespaces}
    
      \AfterEndEnvironment{corollary}{\noindent\ignorespaces}
\def\inj{\mathrm{in}}
\def\sub{\mathrm{sub}}
\newcommand{\coeq}{c}
\def\backgrnd{black!3}	% Background for Tikz pictures

\tikzstyle{place}=
[circle,thick,draw=blue!75,fill=blue!20,minimum size=6mm]
\tikzstyle{transition}=
[rectangle,thick,draw=black!75,fill=black!20,minimum size=4mm]

\newcommand{\Msets}[1]{{#1}^{\oplus}} % Set of multisets on base set #1

\newcommand{\Net}[1]{(\Pl{#1},\Tr{#1},\Pin{-}{#1},\Pout{-}{#1})} % Shorthand for a net
\newcommand{\Pl}[1]{P_{#1}} % Set of places for a net
\newcommand{\Tr}[1]{T_{#1}} % Set of transitions for a net
\newcommand{\Pin}[2]{\mathsf{s}_{#2}(#1)}%{{^\circ}(#1)_{#2}} % Net input of transition #1
\newcommand{\Pout}[2]{\mathsf{t}_{#2}(#1)}%{{(#1)_{#2}^\circ}} % Net output of transition #1
\newcommand{\PlM}[1]{{#1}_P} % Place component of a net morphism
\newcommand{\TrM}[1]{{#1}_T} % Transition component of a net morphism
\newcommand{\GObj}[1]{\operatorname{GenObj} \, #1} % Set of objects of category #1
\newcommand{\GMor}[1]{\operatorname{GenMor} \, #1} % Set of objects of category #1
\newcommand{\Cp}{\fatsemi} % Morphism composition in diagrammatic order

\newcommand{\Source}[1]{\operatorname{s}(#1)} % Domain of function/morphism #1
\newcommand{\Target}[1]{\operatorname{t}(#1)} % Domain of function/morphism #1
\newcommand{\Id}[1]{id_{#1}} % Identity morphism of object #1

\newcommand{\CategoryB}{\mathcal{B}}
\newcommand{\CategoryC}{\mathcal{C}}
\newcommand{\CategoryD}{\mathcal{D}}
\newcommand{\CategoryE}{\mathcal{E}}

\newcommand{\Petri}{\mathbf{Petri}}
\newcommand{\FSSMC}{\mathbf{FSSMC}}
\newcommand{\PetriS}{\Petri^{\Semantics}}
\newcommand{\PetriG}[1]{\Petri^{\Semantics_{#1}}}

\newcommand{\Free}[1]{\mathsf{F}(#1)}
\newcommand{\UnFree}[1]{\mathsf{U}(#1)}

\newcommand{\Fun}[1]{{#1}^\sharp}
\def\FunN{\upsilon_N}
\def\FunM{\upsilon_M}
\def\FunK{\upsilon_K}
\newcommand{\NetSem}[1]{\left(#1, \upsilon_{#1}\right)}

\newcommand{\Semantics}{\mathcal{S}}

\newcommand{\Red}[1]{\mathbf{Red{#1}}}

\newcommand{\Set}{\textbf{Set}} % Category of sets and functions
\newcommand{\Cat}{\textbf{Cat}} % Category of categories

\newcommand{\Tensor}{\otimes} % Monoidal tensor
\newcommand{\Suchthat}[2]{\left\{#1 \: \middle\vert \: #2\right\}} % Set of elements #1 such that condition #2 holds

\tikzset{ % WD-Oriented wiring diagrams
  oriented WD/.style={%everything after equals replaces "oriented WD" in key.
    every to/.style={
      out=0,in=180,draw
    },
    label/.style={
      font=\everymath\expandafter{\the\everymath\scriptstyle},
      inner sep=0pt,
      node distance=2pt and -2pt
    },
    semithick,
    node distance=1 and 1,
    decoration={
      markings, mark=at position \stringdecpos with \stringdec
    },
    ar/.style={
      postaction={decorate}
    },
    execute at begin picture={
      \tikzset{
        x=\bbx, y=\bby,
        every fit/.style={
          inner xsep=\bbx, inner ysep=\bby
        }
      }
    }
  },
  string decoration/.store in=\stringdec,
  string decoration={
    \arrow{stealth};
  },
  string decoration pos/.store in=\stringdecpos,
  string decoration pos=.7,
  bbx/.store in=\bbx,
  bbx = 1.5cm,
  bby/.store in=\bby,
  bby = 1.5ex,
  bb port sep/.store in=\bbportsep,
  bb port sep=1.5,
  bb port length/.store in=\bbportlen,
  bb port length=4pt,
  bb penetrate/.store in=\bbpenetrate,
  bb penetrate=0,
  bb min width/.store in=\bbminwidth,
  bb min width=1cm,
  bb rounded corners/.store in=\bbcorners,
  bb rounded corners=2pt,
  bb small/.style={
    bb port sep=1,
    bb port length=2.5pt,
    bbx=.4cm, bb min width=.4cm,
    bby=.7ex
  },
  bb medium/.style={
    bb port sep=1,
    bb port length=2.5pt,
    bbx=.4cm,
    bb min width=.4cm,
    bby=.9ex
  },
  bb/.code 2 args={%When you see this key, run the code below:
    \pgfmathsetlengthmacro{\bbheight}{\bbportsep * (max(#1,#2)+1) * \bby}
    \pgfkeysalso{
      draw,
      minimum height=\bbheight,
      minimum width=\bbminwidth,
      outer sep=0pt,
      rounded corners=\bbcorners,
      thick,
      prefix after command={
        \pgfextra{\let\fixname\tikzlastnode}
      },
      append after command={
        \pgfextra{
          \draw
          \ifnum #1=0
            {}
          \else
            foreach \i in {1,...,#1} {
              ($(\fixname.north west)!{\i/(#1+1)}!(\fixname.south west)$) +(-
            \bbportlen,0)
            coordinate (\fixname_in\i) -- +(\bbpenetrate,0) coordinate (\fixname_in\i')
            }
          \fi
          \ifnum
            #2=0{}
          \else
            foreach \i in {1,...,#2} {
            ($(\fixname.north east)!{\i/(#2+1)}!(\fixname.south east)$) +(-
            \bbpenetrate,0)
            coordinate (\fixname_out\i') -- +(\bbportlen,0) coordinate (\fixname_out\i)
            }
          \fi;
        }
      }
    }
  },
  bb name/.style={
    append after command={
      \pgfextra{
        \node[anchor=north] at (\fixname.north) {#1}
      ;}
    }
  }
}

\title{The Essence of Petri Net Gluings}

\author{Fabrizio Genovese, Fosco Loregian, and Daniele Palombi}
\address{20[\_]\\[5pt] Tallinn University of Technology\\[5pt] Sapienza University of Rome}
\eaddress{fabrizio.romano.genovese@gmail.com\CR fosco.loregian@taltech.ee\CR danielepalombi@protonmail.com}
\thanks{The second author was supported by the ESF funded Estonian IT Academy research measure (project 2014-2020.4.05.19-0001).}
\copyrightyear{\the\year}

\begin{document}
\maketitle
\begin{abstract}
  Many categorical frameworks have been proposed to formalize the idea of gluing Petri nets with each other. Such frameworks model net gluings in terms of sharing of resources or synchronization of transitions. Interpretations given to these gluings become less satisfactory when we consider Petri nets with a semantics attached to them.

  In this work, we define a framework to compose Petri nets together in such a way that their semantics is respected. In addition to this, we show how our framework generalizes the previously defined ones.
\end{abstract}
\tableofcontents
\section{Introduction}
  \label{sec: introduction}
In the last years, applications of category theory to concurrency,
and in particular to Petri nets, have been flourishing~\cite{Sobocinski2013,
Sobocinski2013a, Fong2016, StateboxTeam2019a,mustcite01}. These
applications span two main directions: On one hand, they organize
Petri nets and their morphisms into categories, and prove interesting
correspondences with other categories~\cite{Meseguer1990, Master2019,9470566}.
Regarding this, much
has been written on how particular classes of symmetric monoidal
categories (SMCs) provide a nice semantics describing the
flow of tokens in a net~\cite{Baldan2003, Sassone1995, Baez2018, Genovese2019}.

On the other hand, researchers tried to employ tools from category
theory to endow Petri nets with ports and connect them together to
form bigger nets. This direction of research can be itself split in multiple
subthreads, which can be categorized as:
\begin{itemize}
  \item Work about connecting nets by merging their
  places~\cite{Baez2018, Fong2016, Baez, Baez2013};
  \item Work about connecting nets by synchronizing their
  transitions~\cite{Bruni2013, Sobocinski2013, Sobocinski2013a};
  \item Other approaches, e.g.~\cite{Baldan2009, Baldan2015,mustcite02}.
\end{itemize}
Perhaps curiously, all these different directions of research
have never been really unified. In particular, mapping nets to
symmetric monoidal categories can be very useful for practical
applications, since it allows to endow Petri nets with a semantics
via monoidal functors~\cite{StateboxTeam2019a}.
Still, it is not really clear what happens
to said semantics when one tries to endow such nets with ports.

Indeed, some net composition paradigms, such
as in~\cite{Bruni2013}, seem to
embrace an observational point of view that
may not be directly compatible with the process interpretation
given by the mapping to SMCs.

On the contrary, place-based net composition as
in~\cite{Fong2016} seems to
be naturally friendlier when it comes to mappings to SMCs.
Still, it surprises that no endeavour, at least to the authors
knowledge, has been made to investigate if the place- and the
transition-based composition paradigms are somehow instances of a more
general notion of net composition.

In this work, we try to answer these questions. By
generalizing the notion of morphism between nets, we show
how different ways of gluing nets together can be modelled
using colimits and functors.
Perhaps surprisingly, we reject the notion of ``ports'' for
nets altogether, arguing that this concept has been somehow
overimposed without really taking net structure into consideration.

Throughout the paper, considerable effort will be made to
point out whether a kind of gluing is or is not
computationally expensive, or ultimately feasible when
it comes to implementing it in code.

The outline of the paper is as follows: In
Section~\ref{sec: Petri nets and free symmetric strict monoidal categories}
we will recap the relationship between Petri nets and free symmetric
monoidal categories on which we are going to rely; in
Section~\ref{sec: the category Petris} we introduce
a new category of Petri nets and morphisms between them;
in Sections~\ref{sec: identifications}
and~\ref{sec: synchronizations} we explain how
place-based and transition-based gluings, respectively,
work in our framework; in
Section~\ref{sec: recovering old concepts gluing between nets}
we will define a monoidal structure and recover the
familiar notions of place- and transition-based gluings
already developed in the literature;
in Section~\ref{sec: changing semantics}
we will show how all the constructions
we defined are parametrized over the choice of a semantics;
in Section~\ref{sec: conclusion and future work}
we conclude by recapping what we did,
defining directions of future work.

\section{Recap: Petri nets and free symmetric strict monoidal categories}
  \label{sec: Petri nets and free symmetric strict monoidal categories}
It is well-known in the category theory folklore
that Petri nets can be considered as free symmetric
strict monoidal categories. Unfortunately, when one tries to
pin down the details of this idea problems arise and different,
inequivalent approaches to associate
a free symmetric strict monoidal category (abbreviated FSSMC)
to a Petri net can be
developed~\cite{Meseguer1990, Sassone1995, Baldan2003, Genovese2019, Baez2018}.
To realize the correspondence, such approaches either use a modified
definition of Petri net~\cite{Baldan2003}, or a
modified definition of FSSMC~\cite{Sassone1995, Baez2018}, or
they weaken the notion of correspondence between the two~\cite{Genovese2019}.

The approach we are going to use has been developed
in~\cite{Genovese2019} and, in contrast to other methods, strives to define
a correspondence between nets and FSSMCs which is computationally
useful and implementable. In~\cite{Genovese2019},
it is extensively discussed
why other approaches can be problematic when implementation
becomes relevant, and we redirect the reader there for more
information (in particular, \cite{Genovese2019} addresses the issue 
of the non-functoriality of the correspondence from Petri nets to FSSMCs).
This is the same approach used also in~\cite{Genovese2021hierar,Genovese2022bdd,Genovese2021mana,Genovese2020guarde}.
% A central assumption to be made for the approach in~\cite{Genovese2019} to
% work is that Petri nets have a well ordering on their
% set of places. This is not a problem since any set can be
% well ordered, while in software applications the places
% of a net are more often than not already enumerated~\cite{StateboxTeam2019}.
% %
% %
% \begin{remark}
%   From now on, we will assume that the set of places
%   of every Petri net is well-ordered.
% \end{remark}
%
We recall the definition of Petri net we'll be working with, from~\cite{Genovese2019}.
\begin{definition}
	A \emph{Petri net} $N$ is a tuple $\Net{N}$, where $\Pl{N}$ and $\Tr{N}$ are sets, called the set of \emph{places} and \emph{transitions} of $N$, respectively, while $\Pin{-}{N}$ and $\Pout{-}{N}$ are functions $\Tr{N} \to \Msets{\Pl{N}}$, representing the input/output places, respectively, connected to each transition.
\end{definition}
\begin{definition}
	A \emph{morphism of Petri nets} $f:N \to M$ is specified by a couple $(\PlM{f}, \TrM{f})$, with $\PlM{f}: \Msets{\Pl{N}} \to \Msets{\Pl{M}}$ multiset homomorphism and $\TrM{f}: \Tr{N} \to \Tr{M}$ function such that
	\begin{equation*}
		\Pin{-}{N};\PlM{f} = \TrM{f};\Pin{-}{M} \qquad \Pout{-}{N};\PlM{f} = \TrM{f};\Pout{-}{M}
	\end{equation*}
\end{definition}
\begin{definition} \label{fac: from Petri to FSSMC}
  To each Petri net $N := \Net{N}$, we associate a
  free strict symmetric monoidal category $\Free{N}$
  defined as follows:
  \begin{itemize}
    \item The places of  $N$ are used to
    freely generate the monoid of objects of $\Free{N}$;
    \item For each transition $t$ in $\Tr{N}$, we define a generating
    morphism $t: \mathfrak{O}(\Pin{t}{N}) \to \mathfrak{O}(\Pout{t}{N})$,
    where $\mathfrak{O}$ denotes the function choosing some ordering on the multiset
    $\Pin{t}{N}$ (and $\Pout{t}{N}$, respectively).% into a string using the
    %order on $\Pl{N}$ to sort elements.
  \end{itemize}
  Details of this construction can be found in~\cite{Genovese2019}.
\end{definition}
Using Fact~\ref{fac: from Petri to FSSMC} we
get a free, symmetric strict
monoidal category $\Free{N}$ for each Petri net $N$.
This FSSMC represents the \emph{category of paths} of the
underlying hypegraph of the net, that is, all the
possible ``allowed'' ways to run it, concurrently, by a
set of actors traversing the hypergraph.
For this reason we call $\Free{N}$ \emph{the category of executions} of $N$.

The correspondence between nets and FSSMCs can be generalized to
a correspondence between the category $\Petri$ of Petri nets and their
morphisms and a suitable category having FSSMCs as its objects:
We are also able to lift an ordinary morphism of Petri nets $M \to N$
to a strict monoidal functor $\Free{M} \to \Free{N}$.
This procedure is, unfortunately, not functorial, meaning that if we
have $f: L \to M$ and $g:M \to N$, then $\Free{f;g}$ may not be equal
(nor isomorphic) to $\Free{f};\Free{g}$ (cf. \cite{Genovese2019} where the issue is addressed). In any case, all functors arising from net
morphisms are identified by the following property:
\begin{remark}%[See \cite{Genovese2019}]
  Given FSSMCs $\CategoryC$ and $\CategoryD$, a functor $F: \CategoryC
  \to \CategoryD$ is called \emph{transition-preserving} if, for each
  morphism generator $t$ of $\CategoryC$, it is $Ft = \sigma\Cp u\Cp \sigma'$,
  with $u$ morphism generator and $\sigma, \sigma'$ symmetries in
  in $\CategoryD$. Every morphism of Petri nets $M \to N$ can be
  lifted to a transition-preserving functor $\Free{M} \to \Free{N}$.
\end{remark}
This was used to prove functoriality from $\FSSMC$ to $\Petri$:
\begin{proposition}%[See \cite{Genovese2019}]
  There is a functor from $\FSSMC$, the
  category having FSSMCs as objects and transition-preserving functors as morphisms, to
  the category $\Petri$ of Petri nets and morphisms between them. This functor
  is denoted with $\UnFree{-}$.
\end{proposition}
\subsection{Folds}
The correspondence between Petri nets and
FSSMCs is particularly useful since, thanks
to the property of being free, FSSMCs can be
easily mapped to other symmetric monoidal categories --
which we interpret as \emph{semantics} -- simply
by specifying maps on the generators.
In the context of our applications, we are only interested in
semantics that can be modelled as symmetric monoidal
categories. We do not deem this to be a restricting
requirement since all process theories, which already
constitute a very broad class, are SMCs.
\begin{remark}
  For the reasons shown above, ``symmetric monoidal category''
  and ``semantics'' will be often treated as interchangeable
  terms in this paper.
  A semantics will be usually denoted as $\Semantics$.
\end{remark}
Unsurprisingly, the mapping to a semantics is realized
via a monoidal functor, which we call a \emph{fold}:
\begin{definition}[See \cite{Genovese2019}]
  Given a net $N$, a \emph{fold}, also called
  an \emph{assignment of semantics for $N$},
  is a strict monoidal functor
  $F: \Free{N} \to \Semantics$,  where $\Semantics$ represents
  a semantics of some sort. Since $\Free{N}$ is free, defining
  where generators are mapped is enough to fully specify $F$.
  Folds are our way to attach meaning to a net.
\end{definition}
For instance, $\Semantics$ can be taken to be the category of
types and functions in some functional programming language:
In this case each generating object is mapped to a
type and each generating morphism to a function. From this,
we obtain a sequence of computations to perform for each
execution of the net.
\section{The category \texorpdfstring{$\PetriS$}{PetriS}}
  \label{sec: the category Petris}
By ``gluing nets together'' we mean
\emph{identifying some places and/or transitions
of a net with the places and/or transitions of some
other net}. Clearly, if an assignment
of semantics is specified for the nets
in question, the gluing should respect it.

Indeed, the idea of identifying places and
transitions differentiates quite a lot when semantics
comes into play, because there are different ways
to reflect the gluing action there. Up to now, this
has been obfuscated by the fact that -- to our
knowledge -- no gluing of nets with assignment
of semantics has been defined or studied in
the literature so far.

Another common point of view in the
literature is the idea of ``gluing places/transitions
of a net with places/transitions of
some \emph{other} net'' -- see e.g.~\cite{Fong2016, Baldan2009}.
This idea seems very natural, but
ultimately obfuscates the principles regulating
net gluings: In short, ``why can't it be that we glue
places and transitions within the \emph{same} net?''
As they are given in their most
basic definition, Petri nets are not endowed with ports
or interfaces of any sort, and the
idea of defining gluings only between
different nets should be considered
as a postulated choice -- which indeed
makes sense in aiding the intuition of
what ``connecting things'' means --
and not as something dictated by the structure.
From this informal consideration it is easy to
realize that many common choices in the literature, such as
endowing a net with left and right ports~\cite{Baez2018,
Baldan2009}, albeit making sense
categorically, do not necessarily
make sense from a net perspective,
where such concepts feel overimposed.

When studying the gluings of nets with semantics in a general setting,
we are somehow forced to move away from these choices, and
embrace the idea that gluing should first be developed naturally
with respect to the net structure: The question should not be
``what do we want to glue?'' but ``what can be glued, and how?''

The first, necessary step to elucidate how
net composition works in a general
setting is to redefine the category $\Petri$
by taking semantics into consideration.
This involves extending the notion of net
morphism by considering generalized
morphisms between their corresponding FSSMCs.
\begin{definition}\label{def: PetriS}
  Fix a semantics $\Semantics$. We define
  the category $\PetriS$ as having:
  \begin{itemize}
    \item As objects, couples $\NetSem{N}$ where
    $N$ is a net and $\FunN: \Free{N} \to \Semantics$ is
    an assignment of semantics;
    \item Morphisms in $\PetriS [\NetSem{M}, \NetSem{N}]$
    consist of all the strict monoidal
    functors $F: \Free{M} \to \Free{N}$ sending object
    generators to object generators,
    and such that $\FunM = F\Cp\FunN$;
    \item Identity on $\NetSem{M}$ is the identity
    functor on $\Free{M}$;
    \item Composition is functor composition.
  \end{itemize}
  Note how we are in no way requiring that
  morphisms in $\PetriS$ are transition-preserving.
  We will often refer to objects of $\PetriS$ just as \emph{nets}.
\end{definition}
It has to be noted that if $\Semantics$
is taken to be the terminal category
then the commutativity condition in
Definition~\ref{def: PetriS} becomes trivial.
An effect of this is that every morphism
in $\Petri$ corresponds,
though not functorially, as explained in~\cite{Genovese2019},
to a morphism in $\PetriS$ via $\Free{-}$.

Simplifying things a bit, permitting
non-transtion-preserving functors
in Definition~\ref{def: PetriS}
has the effect of allowing to map transitions
of a net into \emph{sequences of firings} of
another net. This idea had already been
suggested in~\cite{Meseguer1990}, but has not been further
investigated in recent times, at least
not by the compositionality/applied category
theory crowd. This generalized definition of
net morphism will be used to define net
synchronizations, which we are going to introduce shortly.
On the other hand, since in FSSMCs we interpret generating objects
as places and their monoidal products as markings, the restriction
on $F$ having to send generating objects
to generating objects is justified by the fact
that we did not find any meaningful interpretation of
what it could mean to send a place
of a net to a marking of another.

In the remainder of this document we will
often resort to the well-known
graphical language for Petri nets,
especially to give examples.
But since we are now considering nets
together with their semantic assignments,
the graphical language needs to be enriched
to take this parameter into account.
We do this in Figure~\ref{fig: net with semantics assignment},
where we have decorated each place
and transition of a net with a letter
representing the objects and
morphisms each component of the net, when lifted to a FSSMC generator
using $\Free{-}$, is mapped to.
\begin{figure}[ht!]
  \centering
    \scalebox{0.5}{
\begin{tikzpicture}[node distance=1.3cm,>=stealth',bend angle=45,auto]
  % First net
  %
  \node [transition] (1b) at (-1.5,-1) {$f$};
  \node [place,tokens=0] (2a) at (0,0) {$A$}
    edge [pre] node[swap] {3} (1b);
  \node [place,tokens=0] (2b) at (0,-1)  {$B$}
        edge [pre] (1b);
  \node [place,tokens=0] (2c)  at (0,-2) {$C$}
        edge [pre] node {5} (1b);
  \node [place,tokens=0] (2d) at (0,-3) {$D$};
  \node [transition] (3a) at (1.5,-1)    {$g$}
      edge [pre] node[swap] {2} (2a)
      edge [pre]  (2b)
      edge [pre] node {3} (2c);
  \node [transition] (3b) at (1.5,-3)     {$h$}
      edge [pre] (2c)
      edge [pre] node {4} (2d);
  \node [place,tokens=0] (4a) at (3,-1)  {$E$}
      edge [pre] (3a);
  \node [place,tokens=0] (4b) at (3,-2)  {$F$}
      edge [pre] (3a)
      edge [pre] (3b);
  \node [transition] (5a) at (4.5,-1)  {$k$}
    edge [pre] (4a)
    edge [pre] (4b);
  \begin{pgfonlayer}{background}
    \filldraw [line width=4mm,join=round,\backgrnd]
    (-2,1.5) rectangle (5,-4.5);
  \end{pgfonlayer}
\end{tikzpicture}
}
    \caption{A Petri net decorated with a semantics assignment. The labels $A,\dots,F$ are objects of the chosen semantic monoidal category, while $f,\dots,k$ are morphisms of the appropriate type: for example, $g : A\otimes B\otimes C\to E\otimes F$.}
    \label{fig: net with semantics assignment}
\end{figure}
\subsection{Monoidal structure of \texorpdfstring{$\PetriS$}{PetriS}}
As it will become apparent in the next sections,
to be able to perform operations such as
adding generators to a net one needs a
concept of ``putting nets next to each other'',
that is, a monoidal structure. Evidently,
since we interact at the same time both
with nets and with the free categories arising
from them, we require the monoidal structure
 to cooperate well
with this correspondence. We start by giving
definitions both for nets and for FSSMCs:
\begin{definition}
  \label{def: coproduct of nets}
  Given Petri nets $M, N$, the net $M + N$ is
  defined as follows:
  \begin{itemize}
    \item Places of $M + N$ are the disjoint union of
    places of $M$ and places of $N$, respectively. Concisely:
    \begin{equation*}
      \Pl{M + N} := \Pl{M} \sqcup \Pl{N}
    \end{equation*}
    \item Transitions of $M + N$ are the disjoint union of
    transitions of $M$ and transitions of $N$, respectively. Concisely:
    \begin{equation*}
      \Tr{M + N} := \Tr{M} \sqcup \Tr{N}
    \end{equation*}
    \item Pre-/post-sets are defined in the obvious way:
    \begin{equation*}
      \Pin{t}{M+N} := \begin{cases}
        \Pin{t}{M} &\text{ iff $t \in \Pl{M}$};\\
        \Pin{t}{N} &\text{ iff $t \in \Pl{N}$}.
      \end{cases}
      \qquad
      \Pout{t}{M+N} := \begin{cases}
        \Pout{t}{M} &\text{ iff $t \in \Pl{M}$};\\
        \Pout{t}{N} &\text{ iff $t \in \Pl{N}$}.
      \end{cases}
    \end{equation*}
  \end{itemize}
\end{definition}
\begin{definition}
  \label{def: coproduct of FSSMCs}
  Given FSSMCs $\CategoryC, \CategoryD$, the FSSMC
  $\CategoryC + \CategoryD$ is defined as follows:
  \begin{itemize}
    \item The generating objects of $\CategoryC + \CategoryD$ are
    the disjoint union of the generating objects of $\CategoryC$ and
    $\CategoryD$, respectively. Concisely:
    \begin{equation*}
      \GObj{\CategoryC + \CategoryD} := \GObj{\CategoryC} \sqcup \GObj{\CategoryD}
    \end{equation*}
    \item The generating morphisms of $\CategoryC + \CategoryD$ are
    the disjoint union of the generating morphisms of $\CategoryC$ and
    $\CategoryD$, respectively. Concisely:
    \begin{equation*}
      \GMor{\CategoryC + \CategoryD} := \GMor{\CategoryC} \sqcup \GMor{\CategoryD}
    \end{equation*}
  \end{itemize}
\end{definition}
Definitions~\ref{def: coproduct of nets} and~\ref{def: coproduct of FSSMCs}
interact well with each other, as proven by the following lemma:
\begin{lemma}
  \label{lem: free distributes over coproducts}
 Given nets $M, N$, it is $\Free{M + N} \simeq \Free{M} + \Free{N}$.
\end{lemma}
Moreover, Definitions~\ref{def: coproduct of nets}
and~\ref{def: coproduct of FSSMCs} define, respectively, the coproduct
of nets in $\Petri$ and the product of FSSMCs in $\Cat$. The product in $\Cat$
is also a biproduct in the category of symmetric strict monoidal categories, and
we can use this to define
a monoidal structure on $\PetriS$:
\begin{lemma}\label{lem: monoidal structure of PetriS}
  $\PetriS$ inherits a monoidal structure from the
  biproduct structure in the category of symmetric monoidal categories.
\end{lemma}

\section{Identifications}
  \label{sec: identifications}

\noindent
From now on, the main goal of this paper will
be to characterize morphisms in $\PetriS$ which
can be used to model net gluings.  We start from
identifications. The idea behind identifications
is to conflate places and transitions
\emph{having equal footprints} -- same inputs and outputs.
If two places (or transitions, respectively), get sent
to the same object (resp. morphism) in $\Semantics$,
then they can be identified. Identifications
can thus be interpreted as a way to suppress
redundant information.

The net in Figure~\ref{fig: examples of identification}
exemplifies the concept: In
Figure~\ref{fig: example of identification basic net}.
two places are both mapped to
the object $C$ in $\Semantics$. Since places
are distinguished only by the type of resource
they hold, the fact that both are given type $C$
means that they are holding the same kind of resource,
and as such they can be identified.
This is performed in Figure~\ref{fig: example of identification places}.
Note how both $g$ and $h$ are now
giving and taking, respectively, two tokens to/from $C$:
Merging places should not modify
the overall topology of the net; $g$ used to
output``A resource of type $C$ and a resource
of type $C$", which is the same as ``two resources
of type $C$".

In Figure~\ref{fig: example of identification transitions}
we identify transitions:
Two transitions in
Figure~\ref{fig: example of identification basic net}
are mapped to the same morphism $f$ of $\Semantics$.
From a semantics perspective,
these two transitions do \emph{exactly}
the same thing, and can be conflated into
one. The conflated
transition consumes and produces exactly the same
number of tokens of its conflated components.
Since one can identify both places and
transitions, identifications scale to
entire subnets, as showcased in
Figure~\ref{fig: example of subnet identification}, where the collapsing of subnets is achieved as follows:

\noindent
We now formalize the intuition provided by the
examples into a definition.
\begin{definition}\label{def: identification}
  A Petri net $\NetSem{N}$ is said to be an
  \emph{identification of $\NetSem{M}$} if
  there is a morphism $F: \NetSem{M} \to \NetSem{N}$ such that:
  \begin{itemize}
   \item There is a Petri net $O$, and a couple of
   \emph{transition-preserving} functors
   $l,r: \Free{O} \to \Free{M}$;
   \item $l \Cp \FunM = r \Cp \FunM$;
   \item $F$ is the coequalizer of $l$ and $r$.
  \end{itemize}
  Note that in this case $\FunN$ coincides
  with the arrow arising from the
  universal property of coequalizers. If $\NetSem{N}$
  is an identification for $\NetSem{M}$,
  then we say that $O$, together with $l, r$,
  is a \emph{witness of the identification}.
\end{definition}
In diagrammatic terms, $\NetSem{N}$ is an identification of $\NetSem{M}$ if there exists a diagram of the following shape 
\[\begin{tikzcd}%[>=stealth]
	{\Free{O}} & {\Free{M}} & {\Free{N}} \\
	&& \Semantics
	\arrow["F", from=1-2, to=1-3]
	\arrow["r"', shift right=1, from=1-1, to=1-2]
	\arrow["l", shift left=1, from=1-1, to=1-2]
	\arrow[from=1-2, to=2-3]
	\arrow[from=1-3, to=2-3]
\end{tikzcd}\]
with the property that the diagram $\Free{O} \rightrightarrows \Free{M} \to \Free{N}$ is a coequalizer of categories.
\begin{figure}[!ht]
  \centering
  \begin{subfigure}[t]{0.32\textwidth}\centering
     \scalebox{0.5}{
	\begin{tikzpicture}[node distance=1.3cm,>=stealth',bend angle=45,auto]
		% First net
		%
		\node [place,tokens=0] (1a) at (0,0){$A$};
		\node [transition] (2a) at (1,1)     {$f$}
			edge [pre] (1a);
		\node [transition] (2b) at (1,-1)     {$f$}
			edge [pre] (1a);
		\node [place,tokens=0] (3a) at (2,0)  {$B$}
			edge [pre] (2a)
			edge [pre] (2b);
		\node [transition] (4a) at (3.5,0)  {$g$}
			edge [pre] node {2} (3a);
		\node [place,tokens=0] (5a) at (4.5,1) {$C$}
			edge [pre] (4a);
		\node [place,tokens=0] (5b) at (4.5,-1) {$C$}
			edge [pre] (4a);
		\node [transition] (6a) at (5.5,0)  {$h$}
			edge [pre] (5a)
			edge [pre] (5b);
		\begin{pgfonlayer}{background}
			\filldraw [line width=4mm,join=round,\backgrnd]
      (-0.75,3) rectangle (6.25,-3);
		\end{pgfonlayer}
	\end{tikzpicture}
}
    \caption{A net $N$.}
    \label{fig: example of identification basic net}
  \end{subfigure}
  \begin{subfigure}[t]{0.32\textwidth}\centering
   \scalebox{0.5}{
\begin{tikzpicture}[node distance=1.3cm,>=stealth',bend angle=45,auto]
  % First net
  %
  \node [place,tokens=0] at (0,0)  {$A$};
  \node [transition] (2a) at (1,1)   {$f$}
    edge [pre] (1a);
  \node [transition] (2b) at (1,-1)     {$f$}
    edge [pre] (1a);
  \node [place,tokens=0] (3a) at (2,0)  {$B$}
  edge [pre] (2a)
  edge [pre] (2b);
  \node [transition] (4a) at (3.5,0)  {$g$}
    edge [pre] node {2} (3a);
  \node [place,tokens=0] (5a) at (4.75,0)  {$C$}
    edge [pre] node {2} (4a);
  \node [transition] (6a) at (6,0) {$h$}
    edge [pre] node {2} (5a);
  \begin{pgfonlayer}{background}
    \filldraw [line width=4mm,join=round,\backgrnd]
    (-0.5,3) rectangle (6.5,-3);
  \end{pgfonlayer}
\end{tikzpicture}
}
    \caption{Identification of places.}
    \label{fig: example of identification places}
  \end{subfigure}
 \begin{subfigure}[t]{0.32\textwidth}\centering
   \scalebox{0.5}{
	\begin{tikzpicture}[node distance=1.3cm,>=stealth',bend angle=45,auto]
		% First net
		%
		\node [place,tokens=0] (1a) at (0,0) {$A$};
		\node [transition] (2a) at (1,0)     {$f$}
			edge [pre] (1a);
		\node [place,tokens=0] (3a) at (2,0) {$B$}
			edge [pre] (2a);
		\node [transition] (4a) at (3.5,0)   {$g$}
			edge [pre] node {2} (3a);
		\node [place,tokens=0] (5a) at (4.5,1) {$C$}
			edge [pre] (4a);
		\node [place,tokens=0] (5b) at (4.5,-1) {$C$}
			edge [pre] (4a);
		\node [transition] (6a) at (5.5, 0) {$h$}
			edge [pre] (5a)
			edge [pre] (5b);
		\begin{pgfonlayer}{background}
			\filldraw [line width=4mm,join=round,\backgrnd]
      (-0.75,3) rectangle (6.25,-3);
		\end{pgfonlayer}
	\end{tikzpicture}
}
  \caption{Identification of transitions.}
  \label{fig: example of identification transitions}
\end{subfigure}
  \caption{Examples of identification. In (b) we identify places, using the map coalescing the two copies of $C$ in (a) into the same $C$; in (c) we identify transitions, using the map that coalesces the two copies of $f$ in (a) into the same $f$.}
  \label{fig: examples of identification}
\end{figure}
\begin{figure}[!ht]
  \centering
  \begin{subfigure}[t]{0.32\textwidth}\centering
    \scalebox{0.5}{
  \begin{tikzpicture}[node distance=1.2cm,>=stealth',bend angle=45,auto]
    \draw[draw=red!5,fill=red!10] (0.5,0.5) rectangle (4.5,3);
    \draw[draw=red!5,fill=red!10] (0.5,-0.5) rectangle (4.5,-3);   
		\node [place,tokens=0] at (0,0) {$A$};
    \node [transition] (2a) at (1,1)    {$f$}
      edge [pre] (1a);
    \node [transition] (2b) at (1,-1)   {$f$}
      edge [pre] (1a);
		\node [place,tokens=0] (3b) at (2.5,1) {$C$}
      edge [pre] node {2} (2a);
    \node [place,tokens=0] (3a) at (2.5,2.5) {$B$}
      edge [post] (2a);
    \node [place,tokens=0] (3c) at (2.5,-1)  {$C$}
      edge [pre] node[swap] {2} (2b);
    \node [place,tokens=0] (3d) at (2.5,-2.5)  {$B$}
      edge [post] (2b);
		\node [transition] (4a) at (4,1)  {$g$}
      edge [pre] node {2} (3b)
      edge [post] (3a);
    \node [transition] (4b) at (4,-1)  {$g$}
      edge [pre] node[swap] {2} (3c)
      edge [post] (3d);
		\node [place,tokens=0] (5a) at (5,0)  {$D$}
      edge [pre] (4a)
      edge [pre] (4b);
		\node [transition] (6a) at (6,0) {$h$}
			edge [pre] (5a);
		\begin{pgfonlayer}{background}
			\filldraw [line width=4mm,join=round,\backgrnd]
      (-0.5,3) rectangle (6.5,-3);
		\end{pgfonlayer}
  \end{tikzpicture}
}
    \caption{Subnets highlighted.}
    \label{fig: example of subnet identification basic net}
  \end{subfigure}
  \begin{subfigure}[t]{0.32\textwidth}\centering
    \scalebox{0.5}{
  \begin{tikzpicture}[node distance=1.2cm,>=stealth',bend angle=45,auto]
    \draw[draw=red!5,fill=red!10] (0.5,0.5) rectangle (4.5,3);
    \draw[draw=red!5,fill=red!10] (0.5,-0.5) rectangle (4.5,-3);
		\node [place,tokens=0] at (0,0) {$A$};
    \node [transition, fill=green!15] (2a) at (1,1)    {$f$}
      edge [pre] (1a);
    \node [transition, fill=green!15] (2b) at (1,-1)   {$f$}
      edge [pre] (1a);
		\node [place,tokens=0, fill=yellow!35] (3b) at (2.5,1) {$C$}
      edge [pre] node {2} (2a);
    \node [place,tokens=0, fill=yellow!15] (3a) at (2.5,2.5) {$B$}
      edge [post] (2a);
    \node [place,tokens=0, fill=yellow!35] (3c) at (2.5,-1)  {$C$}
      edge [pre] node[swap] {2} (2b);
    \node [place,tokens=0, fill=yellow!15] (3d) at (2.5,-2.5)  {$B$}
      edge [post] (2b);
		\node [transition, fill=green!30] (4a) at (4,1)  {$g$}
      edge [pre] node {2} (3b)
      edge [post] (3a);
    \node [transition, fill=green!30] (4b) at (4,-1)  {$g$}
      edge [pre] node[swap] {2} (3c)
      edge [post] (3d);
		\node [place,tokens=0] (5a) at (5,0)  {$D$}
      edge [pre] (4a)
      edge [pre] (4b);
		\node [transition] (6a) at (6,0) {$h$}
			edge [pre] (5a);
		\begin{pgfonlayer}{background}
			\filldraw [line width=4mm,join=round,\backgrnd]
      (-0.5,3) rectangle (6.5,-3);
		\end{pgfonlayer}
  \end{tikzpicture}
}
    \caption{Highlighting components.}
    \label{fig: example of subnet identification middle1 net}
  \end{subfigure}
  \begin{subfigure}[t]{0.32\textwidth}\centering
   \scalebox{0.5}{
  \begin{tikzpicture}[node distance=1.2cm,>=stealth',bend angle=45,auto]
    \draw[draw=red!5,fill=red!10] (0.5,-2.5) rectangle (4.5,2.5);
		\node [place,tokens=0] at (0,0) {$A$};
    \node [transition] (2a) at (1,0)   {$f$}
      edge [pre] (1a);
		\node [place,tokens=0] (3b) at (2.5, -1.5)  {$C$}
      edge [pre] node {2} (2a);
    \node [place,tokens=0] (3a) at (2.5,1.5)  {$B$}
      edge [post] (2a);
		\node [transition] (4a) at (4,0)  {$g$}
      edge [pre] node {2} (3b)
      edge [post] (3a);
		\node [place,tokens=0] (5a) at (5,0) {$D$}
      edge [pre] (4a)	;
		\node [transition] (6a) at (6,0)  {$h$}
			edge [pre] (5a);
		\begin{pgfonlayer}{background}
			\filldraw [line width=4mm,join=round,\backgrnd]
      (-0.5,3) rectangle (6.5,-3);
		\end{pgfonlayer}
  \end{tikzpicture}
}
    \caption{Subnets identified.}
    \label{fig: example of subnet identification after}
  \end{subfigure}
  \caption{Identification of subnets. Here a whole subnet (highlighted in (a)) is identified in two steps; first, identifiying places, and then identifying transitions. The places and transitions identified are highlighted in matching colors in (b). The result is shown in (c).}
  \label{fig: example of subnet identification}
\end{figure}
In this definition, the net $O$ is used
to select the places and transitions to be identified
in the target net $M$. The functors $l,r$
are required to be transition-preserving because
components of $O$ should not be mapped to
computations of $M$: The mapping acts
at a topological level, literally identifying places
and transitions of the net, not their executions. This amounts
to ask that $l,r$ arise from maps between
the underlying hypergraphs of $O$ and $M$.

On the other hand, the condition $l \Cp \FunM = r \Cp \FunM$
formalizes the requirement that the places and
transitions to be identified correspond to the
same objects and morphisms in the semantics $\Semantics$.
This is conceptually obvious, as,
in our framework, it does not make sense to identify
components mapped to different semantic entities.

Having used $O$ to select the components of $M$ to be
identified, the actual identification is performed
by using coequalizers. The fact that $l,r$ are
transition-preserving provides a
constructive characterization:
\begin{lemma}\label{lem: characterizing coequalizers}
  Let $\CategoryC, \CategoryD$ be free strict
  symmetric monoidal categories, and let
  $F,G: \CategoryC \to \CategoryD$ be a couple
  of transition-preserving functors
  sending generating objects to generating objects.
  Denoting with $\GObj{\CategoryC}$
  and $\GMor{\CategoryC}$ the generating
  objects and morphisms, respectively, of $\CategoryC$, then the
  coequalizer $\CategoryE$ of $F,G$ is the
  following free symmetric strict monoidal category:
  \begin{itemize}
    \item Generating objects of $\CategoryE$ are
    $\GObj{\CategoryD}/\simeq_o$, where $\simeq_o$
    is the equivalence relation generated by
    \begin{equation*}
      \forall x \in \GObj{\CategoryC}.(Fx = Gx)
    \end{equation*}
    \item Given a generating morphism $f$ of $\CategoryC$, and
    denoting with $f_F$ and $f_G$ the generating morphisms
    of $\CategoryD$ such that $Ff = \sigma \Cp f_F \Cp \sigma'$ and
   $Gf = \varsigma \Cp f_G \Cp \varsigma'$, generating
   morphisms of $\CategoryE$ are
   $\GMor{\CategoryD}/\simeq_m$, where $\simeq_m$
   is the equivalence relation generated by
   \begin{equation*}
    \forall f \in \GMor{\CategoryC}.(f_F = f_G)
   \end{equation*}
 \end{itemize}
\end{lemma}
\subsection{Implementation perspective}
Lemma~\ref{lem: characterizing coequalizers}
provides a constructive way to build coequalizers,
but it is, unfortunately, still computationally unfeasible.
It makes heavy use of equivalence relations,
and quotient-like structures are notoriously
problematic to implement~\cite{Licata2011}.
Luckily, the most useful way to employ
identifications is to merge a finite number of
places together, and in this setting a
computationally friendly characterization of
the coequalizers involved can be given.

Indeed, it has to be noted that the
Petri net consisting of just one
place and no transitions can be used as a witness
to merge two places of a net together.
The advantage of this approach is that
calculating the coequalizer becomes
a very easy task:
\begin{lemma}\label{lem: identifying two places}
 Let $O_1$ be the Petri net having only one
 place, $o$, and no transitions.
 Let furthermore $\CategoryC$ be a free symmetric
 strict monoidal category, and let
 $F,G: \Free{O_1} \to \CategoryC$ be a couple of
 strict monoidal functors between them.
 Then the coequalizer $\CategoryD$ of $F,G$ is the
 following FSSMC:
 \begin{itemize}
   \item Generating objects of $\CategoryD$ are $\GObj{\CategoryC} \backslash \{Go\}$, where $\backslash$ denotes
   subtraction of sets;
   \item Generating morphisms of $\CategoryD$ are:
   \begin{equation*}
    \{\alpha': \Source{\alpha}[Go \backslash Fo] \to \Target{\alpha}[Go \backslash Fo] \mid \alpha \in \GMor{\CategoryC} \}
   \end{equation*}
   That is, the generating morphisms of $\CategoryD$ are the generating
   morphisms of $\CategoryC$ where any occurrence of $Go$ in their
   source and target is substituted with $Fo$.
 \end{itemize}
 Moreover, the coequalizing morphism $\CategoryC \to \CategoryD$ is
 transition-preserving.
\end{lemma}
Lemma~\ref{lem: identifying two places} means
that when identifying two places $A$ and $B$ together,
it is enough to strip one of the
generating objects from $\CategoryC$ and rename
all of its occurrencies with the other.
The following result shows how
Lemma~\ref{lem: identifying two places} can be
scaled to any finite number of places:
\begin{lemma}\label{lem: identifying n places}
 Denote with $O_n$ the Petri net consisting of $n$ places
 and no transitions.
 Let $\NetSem{N}$ be an identification of $\NetSem{M}$ via
 $F: \NetSem{M} \to \NetSem{N}$ with witness $O_n,l,r$.

 Then there exist transition-preserving functors
 $l_1, r_1, \dots, l_n, r_n$ such that the following
 diagram commutes, where $\coeq_i$ is the coequalizer
 of $l_i, r_i$, and dashed arrows are obtained
 from the universal property of coequalizers:
 \begin{equation*}
   % \scalebox{0.8}{
    \begin{tikzpicture}[node distance=2cm,>=stealth',bend angle=45,auto,xscale=2.5cm]
        \node (1) at (0,0) {$\Free{M}$};
        \node (2) [right of=1] {$\Free{N}_1$};
        \node (3) [right of=2] {$\dots$};
        \node (n-1) [right= 1.5cm of 3] {$\Free{N}_{n-1}$};
        \node (n) [right= 1.5cm of n-1] {$\Free{N}$};
        \node (w) [left of=1] {$\Free{O_n}$};

        \node(1a) [above of=1] {$\Free{O_1}$};
        \node(2a) [above of=2] {$\Free{O_1}$};
        \node(na) [above of=n-1] {$\Free{O_1}$};

        \node(1b) [below of=1] {$\Semantics$};

        \draw[transform canvas={yshift=0.5ex},->] (w) to node[font=\tiny] {$l$} (1);
        \draw[transform canvas={yshift=-0.5ex},->](w) to node[font=\tiny,swap] {$r$} (1);
        \draw[transform canvas={xshift=0.5ex},->] (1a) to node[font=\tiny] {$l_1$} (1);
        \draw[transform canvas={xshift=-0.5ex},->](1a) to node[font=\tiny,swap] {$r_1$} (1);
        \draw[transform canvas={xshift=0.5ex},->] (2a) to node[font=\tiny] {$l_2$} (2);
        \draw[transform canvas={xshift=-0.5ex},->](2a) to node[font=\tiny,swap] {$r_2$} (2);
        \draw[transform canvas={xshift=0.5ex},->] (na) to node[font=\tiny] {$l_n$} (n-1);
        \draw[transform canvas={xshift=-0.5ex},->](na) to node[font=\tiny,swap] {$r_n$} (n-1);

        \draw[dotted, ->] (1) to node[font=\tiny] {$\coeq_1$} (2);
        \draw[dotted, ->] (2) to node[font=\tiny] {$\coeq_2$} (3);
        \draw[dotted, ->] (3) to node[font=\tiny] {$\coeq_{n-1}$} (n-1);
        \draw[dotted, ->] (n-1) to node[font=\tiny] {$\coeq_{n}$} (n);

        \draw[->] (1) to node[font=\tiny,swap] {$\FunM$} (1b);
        \draw[dashed, ->] (2) to node[font=\tiny] {} (1b);
        \draw[dashed, ->] (n-1) to node[font=\tiny] {} (1b);

        \draw[->] (n) to node[font=\tiny] {$\FunN$} (1b);
        \end{tikzpicture}
% }
 \end{equation*}
 Moreover, it is $F = \coeq_1 \Cp \dots \Cp \coeq_{n}$. In other
 words: The identification of a finite,
 arbitrary number of places can be
 performed in steps, where at each step
 no more than two places are identified.
\end{lemma}
These two lemmas give a very covenient way
of identifying places, which does not
involve resorting to quotients. Moreover,
identifications can be constructed easily
just by specifying witnesses.
Consider, in fact, a net $\NetSem{M}$:
To construct an identification for it,
it is sufficient to specify a net $W$ and functors
$l,r: \Free{W} \to \Free{M}$. Notably,
no requirement on the semantics of $W$ is made,
and $l, r$ are considered ``acceptable''
if they satisfy the condition $l \Cp \FunM = r \Cp \FunM$.
When $W = O_n$ for some $n$,
Lemma~\ref{lem: identifying n places}
can be applied to
split the computation in steps, and
the coequalizer at each step can be calculated
using Lemma~\ref{lem: identifying two places}.
Composing the results, one gets the coequalizer of $l, r$,
namely a functor
$F:\Free{M} \to \CategoryC$, for
some free symmetric strict monoidal
category $\CategoryC$. $\UnFree{\CategoryC}$
will then be the resulting net, and $\Fun{\UnFree{\CategoryC}}$
can be obtained automatically considering that, since
$l \Cp \FunM = r \Cp \FunM$, the universal
property of coequalizers implies the existence
of a unique arrow
$\Fun{\UnFree{\CategoryC}}: \Free{\UnFree{\CategoryC}} \simeq\CategoryC \to \Semantics$ such
that $F \Cp  \Fun{\UnFree{\CategoryC}} = \FunM$,
instantly turning $F$ into a morphism
$\NetSem{M} \to \NetSem{\UnFree{\CategoryC}}$.
We conclude that $\NetSem{M}$ and witnesses
$W,l,r$ are enough to fully specify the
identification $\NetSem{\UnFree{\CategoryC}}$.
\section{Synchronizations}
  \label{sec: synchronizations}
Now we focus on another interesting class of
morphisms in $\PetriS$, \emph{synchronizations}.
The underlying idea of synchronizations is
that we would like to compress multiple
events happening in a net to one: For instance,
we would like to say that each time a
transition fires in a net, this automatically
triggers the firing of some other transition,
which could be either located in a different
net or not.
In short, this calls for multiple, separated
firings to be conflated into a unique event.
\begin{figure}[!ht]
  \centering
  \begin{subfigure}[t]{0.49\textwidth}\centering
    \scalebox{0.5}{
\begin{tikzpicture}[node distance=1.3cm,>=stealth',bend angle=45,auto]
  % First net
  %
  \node [transition] (1b) at (-1.5,-1) {$f$};
  \node [place,tokens=0] (2a) at (0,0) {$A$}
    edge [pre] node[swap] {3} (1b);
  \node [place,tokens=0] (2b) at (0,-1)  {$B$}
        edge [pre] (1b);
  \node [place,tokens=0] (2c)  at (0,-2) {$C$}
        edge [pre] node {5} (1b);
  \node [place,tokens=0] (2d) at (0,-3) {$D$};
  \node [transition] (3a) at (1.5,-1)    {$g$}
      edge [pre] node[swap] {2} (2a)
      edge [pre]  (2b)
      edge [pre] node {3} (2c);
  \node [transition] (3b) at (1.5,-3)     {$h$}
      edge [pre] (2c)
      edge [pre] node {4} (2d);
  \node [place,tokens=0] (4a) at (3,-1)  {$E$}
      edge [pre] (3a);
  \node [place,tokens=0] (4b) at (3,-2)  {$F$}
      edge [pre] (3a)
      edge [pre] (3b);
  \node [transition] (5a) at (4.5,-1)  {$k$}
    edge [pre] (4a)
    edge [pre] (4b);
  \begin{pgfonlayer}{background}
    \filldraw [line width=4mm,join=round,\backgrnd]
    (-2,1.5) rectangle (5,-4.5);
  \end{pgfonlayer}
\end{tikzpicture}
}
    \caption{A net.}
    \label{fig: example of synchronization basic net}
 \end{subfigure}
 \begin{subfigure}[t]{0.49\textwidth}\centering
   \scalebox{0.5}{
\begin{tikzpicture}[node distance=1.3cm,>=stealth',bend angle=45,auto]
  % First net
  %
  \node [transition] (1b) at (-1.5,-1) {$f$};
  \node [place,tokens=0] (2a) at (0,0) {$A$}
    edge [pre] node[swap] {3} (1b);
  \node [place,tokens=0] (2b) at (0,-1)  {$B$}
        edge [pre] (1b);
  \node [place,tokens=0] (2c)  at (0,-2) {$C$}
        edge [pre] node {5} (1b);
  \node [place,tokens=0] (2d) at (0,-3) {$D$};
  \node [transition] (3a) at (1.5,-1)   {$g \Cp k$}
      edge [pre] node[swap] {2} (2a)
      edge [pre]  (2b)
      edge [pre] node {3} (2c);
  \node [transition] (3b) at (1.5,-3)   {$h$}
      edge [pre] (2c)
      edge [pre] node {4} (2d);
  \node [place,tokens=0] (4a) at (3,-1)  {$E$};
  \node [place,tokens=0] (4b) at (3,-2)  {$F$}
      edge [pre] (3b);
  \begin{pgfonlayer}{background}
    \filldraw [line width=4mm,join=round,\backgrnd]
    (-2.625,1.5) rectangle (4.325,-4.5);
  \end{pgfonlayer}
\end{tikzpicture}
}
   \caption{Synchronization of $g,k$, resources shared.}
   \label{fig: example of synchronization resources shared}
 \end{subfigure}

 \bigskip
 \begin{subfigure}[t]{0.49\textwidth}\centering
   \scalebox{0.5}{
\begin{tikzpicture}[node distance=1.3cm,>=stealth',bend angle=45,auto]
  % First net
  %
  \node [transition] (1b) at (-1.5,-1) {$f$};
  \node [place,tokens=0] (2a) at (0,0) {$A$}
    edge [pre] node[swap] {3} (1b);
  \node [place,tokens=0] (2b) at (0,-1)  {$B$}
        edge [pre] (1b);
  \node [place,tokens=0] (2c)  at (0,-2) {$C$}
        edge [pre] node {5} (1b);
  \node [place,tokens=0] (2d) at (0,-3) {$D$};
  \node [transition] (3a) at (1.5,-2)    {$g \Tensor h$}
      edge [pre] node[swap] {2} (2a)
      edge [pre]  (2b)
      edge [pre] node {4} (2c)
      edge [pre] (2c)
      edge [pre] node {4} (2d);
  \node [place,tokens=0] (4a) at (3,-1)  {$E$}
      edge [pre] (3a);
  \node [place,tokens=0] (4b) at (3,-2)  {$F$}
      edge [pre] node {$2$} (3a);
  \node [transition] (5a) at (4.5,-1)  {$k$}
    edge [pre] (4a)
    edge [pre] (4b);
  \begin{pgfonlayer}{background}
    \filldraw [line width=4mm,join=round,\backgrnd]
    (-2,1.5) rectangle (5,-4.5);
  \end{pgfonlayer}
\end{tikzpicture}
}
   \caption{Synchronization of $g,h$, resources not shared.}
   \label{fig: example of synchronization resources not shared}
 \end{subfigure}
 \begin{subfigure}[t]{0.49\textwidth}\centering
    \scalebox{0.5}{
\begin{tikzpicture}[node distance=1.3cm,>=stealth',bend angle=45,auto]
  % First net
  %
  \node [transition] (1b) at (-1.5,-1) {$f$};
  \node [place,tokens=0] (2a) at (0,0) {$A$}
    edge [pre] node[swap] {3} (1b);
  \node [place,tokens=0] (2b) at (0,-1)  {$B$}
        edge [pre] (1b);
  \node [place,tokens=0] (2c)  at (0,-2) {$C$}
        edge [pre] node {5} (1b);
  \node [place,tokens=0] (2d) at (0,-3) {$D$};
  \node [transition] (3a) at (3,-2)  {$(g \Tensor h) \Cp (k \Tensor \Id{F})$}
      edge [pre] node[swap] {2} (2a)
      edge [pre]  (2b)
      edge [pre] node {4} (2c)
      edge [pre] (2c)
      edge [pre] node {4} (2d);
  \node [place,tokens=0] (4a) at (3,-0.5)  {$E$};
  \node [place,tokens=0] (4b) at (3,-3.5) {$F$}
      edge [pre] (3a);
  \begin{pgfonlayer}{background}
    \filldraw [line width=4mm,join=round,\backgrnd]
    (-2,1.5) rectangle (5,-4.5);
  \end{pgfonlayer}
\end{tikzpicture}
}
    \caption{Synchronization of $g,k,h$.}
    \label{fig: example of synchronization mixed}
 \end{subfigure}
 \caption{Examples of net synchronizations.}
  \label{fig: example of synchronization}
\end{figure}

\noindent
When considering the process semantics
given by FSSMCs, synchronizations can happen
in at least two ways: They can share resources
or not. We elucidate the differences between
the two approaches with examples, starting
with sharing of resources. In
Figure~\ref{fig: example of synchronization resources shared},
two transitions of the net in
Figure~\ref{fig: example of synchronization basic net}
have been conflated.
The conflated firing sequence prescribes
that one transition immediately consumes the
resources produced by the other. In this
case, if the to-be-synced transitions get mapped to
morphisms $g$ and $k$ in $\Semantics$,
respectively, then it is clear that the conflated
transition has to be mapped to $g \Cp k$.
As it is evident, the place $E$ is now
isolated; this is in line with the
intuition that the resources produced by
$g$ are \emph{instantly consumed}
by $h$ because the firings are happening
synchronously, and hence do not pass through $E$.
Moreover, from the example one infers how
this kind of synchronization is not
always possible since some obvious
compatibility conditions on the transitions'
domains and codomains
have to be required.

The other kind of synchronization one
can consider does not involve
sharing of resources, merely
using one transition to trigger the
firing of another. Interestingly, since
synchronizations between firings
are thought of
as instantaneous, the information about
which transition triggers the firing
can be discarded altogether.
This is shown in
Figure~\ref{fig: example of synchronization resources not shared},
and in this case it is clear that
if the to-be-synced transitions
are mapped to $g$, $h$ in $\Semantics$,
respectively, then the conflated
transition has to be mapped to
their monoidal product.
Since no sharing of resources is involved, this
kind of gluing has no preconditions and any
number of transitions in the same net can be
synchronized.

Finally, these two kinds of
synchronization can be combined, as it is shown in
Figure~\ref{fig: example of synchronization mixed}.
Notice how some extra bureaucracy
is needed here:
The synchronization $g \Tensor h$
produces two tokens in the place labelled $F$,
while the transition labelled with $k$ will
consume just one. This means that after firing
there will be one token in $F$ -- the one that $k$
did not consume -- and it will have to be
explicitly declared which token of the two this is.

Having provided intuition of how synchronizations
work, it is useful to draw comparisons with
identifications, defined in the previous
section. Identifications are better understood
as a way to eliminate redundancies and, in
particular, when applied to places they model
``asynchronous'' gluings;
on the contrary, synchronizations
act at a ``message level'', coordinating the firings
of different transitions.
Whereas identifications
both apply to transitions and places -- since both
can be redundant -- synchronizations
operate solely at the transition level:
As we already stressed places, thought of as
as ``bags'' holding resources, are
inherently asynchronous, making the concept
of synchronization meaningless in this context.

We now want to define synchronizations formally.
These can be understood
as multiple things happening at the same time, which we will model separately:
\begin{itemize}
  \item Generators corresponding to the
  sychronized transitions have to be added to the net;
  \item Generators corresponding to the
  transitions defining the synchronization must be erased.
\end{itemize}
These two steps will be formalized using pushout squares,
the final characterization of synchronization being
very similar to the application of a
double pushout rewriting rule~\cite{habel_muller_plump_2001}.

As for identifications, a witness is needed
to perform the synchronization. This witness
will ultimately determine which generators have
to be erased and which computations have to be promoted
to generators, or, in other words, will completely
characterize the double pushout rewriting rule to be applied.
\begin{definition}
  A \emph{synchronization
  witness for $\NetSem{N}$} is a pair $(W,w)$ where $W$ is a Petri
  net $W$ and $w$ a strict monoidal functor
  $\Free{W} \to \Free{N}$ which sends
  generating objects to generating objects,
  is injective on objects and faithful.
\end{definition}
Note how, in contrast with identifications,
we do not require in any way $w$ to
be transition-preserving, meaning that
morphism generators of $\Free{W}$ may not
be morphism generators of $\Free{N}$.
Another way to see this is that
even if $\Free{W}$ can be thought of as a
subcategory of $\Free{N}$,
$W$ may not be a subnet of $N$.

The following definition will also be used
heavily in the remainder of this section,
and hence deserves a special notation:
\begin{definition}
  Given a net $W$, we denote with $\overline{W}$ the net
  having the same places of $W$ but no transitions. There
  is an obvious identity on objects functor $\inj_W:\Free{\overline{W}}
  \hookrightarrow \Free{W}$, called \emph{inclusion
  of $\Free{\overline{W}}$ into $\Free{W}$}.
\end{definition}
\subsection{Promoting computations to generators}
In implementing synchronizations, it is
necessary to promote entire computations of
a net to being morphism generators.
Such newly added generators
represent the act of conflating events as intuitively
described in Figure~\ref{fig: example of synchronization}.
\begin{definition}
  \label{def: addition of generators}
  Given a net $\NetSem{K}$ and a synchronization
  witness $W,w$ for $\NetSem{K}$, a net $\NetSem{M}$ is an
  \emph{addition of generating morphisms to $\NetSem{K}$
  via $W, w$}
  if $\Free{M}$ is the pushout of
  $\Free{\overline{W}} \xhookrightarrow{\inj_W} \Free{W} \xrightarrow{w} \Free{K}$
  and $\Free{\overline{W}} \xhookrightarrow{\inj_W} \Free{W}$, and $\FunM$ arises from
  the universal property of the pushout.
\end{definition}
Definition~\ref{def: addition of generators} is rather succint.
To better understand why it models the idea of adding generators
to a net, we notice that the pushout can be also characterized as
the coequalizer of a coproduct.
Indeed, given a net $\NetSem{K}$  and a synchronization
witness $W, w:\Free{W} \to \Free{K}$, the
coproduct structure can be used to join $W$ to $K$,
obtaining the net $(K + W, [\FunK, w \Cp \FunK])$:
\begin{equation*}
  \scalebox{0.75}{
  \begin{tikzpicture}[node distance=4cm,>=stealth',bend angle=45,auto]
    \node (1) at (0,0) {$\Free{K + W}$};
    \node (2) [right of=1] {};
    \node (3) [right of=2] {$\Free{K}$};
    \node (4) [below of=1] {$\Free{W}$};
    \node (5) [below of=2] {$\Free{K}$};
    \node (6) [below of=3] {$\Semantics$};
    \draw[left hook-latex] (3) to node[swap] {$\iota_1$} (1);
    \draw[right hook-latex] (4) to node {$\iota_2$} (1);
    \draw[->] (4) to node[swap] {$w$} (5);
    \draw[->] (3) to node {$\FunK$} (6);
    \draw[->] (5) to node[swap] {$\FunK$} (6);
    \draw[dashed, ->] (1) to node[swap] {$[\FunK, w \Cp \FunK]$} (6);
  \end{tikzpicture}
}
\end{equation*}
In doing so, along with adding the
morphisms of $W$ to $K$, we also took
the disjoint union of their objects, which are
now redundant. To understand this, look at
Figure~\ref{fig: adding and merging generators to a net}:
Starting from the net in
Figure~\ref{fig: example of synchronization starting net},
we defined a witness $W, w$ pinpointing the computation
$f \Cp g$, and then promoted such computation
from being just a morphism in $\Free{K}$ to being a
morphism generator by taking the coproduct $\Free{K +W}$,
as it is shown in
Figure~\ref{fig: example of synchronization added net}.
This generator though ``floats
on its own'', since the generating objects of $\Free{W}$
have been added to $K$ as well.
Even if they are mapped by $w$ to
generating objects of $K$, the coproduct
structure considers them as separated and
makes them redundant.
To obtain the wanted result depicted in
Figure~\ref{fig: example of synchronization merged net}, we need
to merge the redundant places introduced by the coproduct.
\begin{figure}[!ht]
  \centering
  \begin{subfigure}[t]{0.32\textwidth}\centering
    \scalebox{0.5}{
  \begin{tikzpicture}[node distance=1.3cm,>=stealth',bend angle=45,auto]
		% First net
		%
		\node [place] (1a) at (0,-1.5) {$A$};
		\node [transition] (2a) at (1.5,1) {$f$}
	  	edge [pre] (1a);
		\node [place] (3a) at (3,1) {$B$}
			edge [pre] (2a);
    \node [transition] (3b) at (3,-4)  {$h$}
     edge [pre] (1a);
    \node [transition] (4a) at (4.5,1)     {$g$}
      edge [pre] (3a);
    \node [place] (5a) at (6,1)  {$C$}
     edge [pre] (4a);
    \node [place] (5b) at (6,-4)    {$D$}
      edge [pre] (3b);
		\begin{pgfonlayer}{background}
			\filldraw [line width=4mm,join=round,\backgrnd]
				(-0.5,1.5) rectangle (6.5,-4.5);
		\end{pgfonlayer}
  \end{tikzpicture}
  }
    \caption{A net $K$.}
    \label{fig: example of synchronization starting net}
  \end{subfigure}
  \begin{subfigure}[t]{0.32\textwidth}\centering
    \scalebox{0.5}{
  \begin{tikzpicture}[node distance=1.3cm,>=stealth',bend angle=45,auto]
		% First net
		%
		\node [place] (1a) at (0,-1.5) {$A$};
    \node [place] (1b) at (1.25,-1.5) {$A$};
		\node [transition] (2a) at (1.5,1) {$f$}
	  	edge [pre] (1a);
		\node [place] (3a) at (3,1) {$B$}
			edge [pre] (2a);
    \node [transition] (3b) at (3,-4)  {$h$}
     edge [pre] (1a);
    \node [transition] (3c) at (3,-1.5)     {$f \Cp g$}
      edge [pre] (1b);
    \node [transition] (4a) at (4.5,1)     {$g$}
      edge [pre] (3a);
    \node [place] (5a) at (6,1)  {$C$}
     edge [pre] (4a);
    \node [place] (5c) at (5,0)  {$C$}
     edge [pre] (3c);
    \node [place] (5b) at (6,-4)    {$D$}
      edge [pre] (3b);
		\begin{pgfonlayer}{background}
			\filldraw [line width=4mm,join=round,\backgrnd]
				(-0.5,1.5) rectangle (6.5,-4.5);
		\end{pgfonlayer}
  \end{tikzpicture}
  }
    \caption{Computation promoted.}
    \label{fig: example of synchronization added net}
  \end{subfigure}
  \begin{subfigure}[t]{0.32\textwidth}\centering
    \scalebox{0.5}{
  \begin{tikzpicture}[node distance=1.3cm,>=stealth',bend angle=45,auto]
		% First net
		%
		\node [place] (1a) at (0,-1.5) {$A$};
		\node [transition] (2a) at (1.5,1) {$f$}
	  	edge [pre] (1a);
		\node [place] (3a) at (3,1) {$B$}
			edge [pre] (2a);
    \node [transition] (3b) at (3,-4)  {$h$}
     edge [pre] (1a);
    \node [transition] (3c) at (3,-1.5)     {$f \Cp g$}
      edge [pre] (1a);
    \node [transition] (4a) at (4.5,1)     {$g$}
      edge [pre] (3a);
    \node [place] (5a) at (6,1)  {$C$}
     edge [pre] (4a)
     edge [pre] (3c);
    \node [place] (5b) at (6,-4)    {$D$}
      edge [pre] (3b);
		\begin{pgfonlayer}{background}
			\filldraw [line width=4mm,join=round,\backgrnd]
				(-0.5,1.5) rectangle (6.5,-4.5);
		\end{pgfonlayer}
  \end{tikzpicture}
  }
    \caption{Generator merged.}
    \label{fig: example of synchronization merged net}
  \end{subfigure}
  \caption{Adding and merging generators to a net.}
    \label{fig: adding and merging generators to a net}
\end{figure}

\noindent
We do
so by using the tools developed in Section~\ref{sec: identifications},
specifically using $\overline{W}$ and the functors $\inj_W \Cp w \Cp \iota_1$, $\inj_W \Cp \iota_2$
as witnesses to perform an identification.
\begin{itemize}
  \item The functors $\Free{\overline{W}} \xhookrightarrow{\inj_W} \Free{W}
  \xrightarrow{w} \Free{K} \xhookrightarrow{\iota_1} \Free{K +W}$ and
  $\Free{\overline{W}} \xhookrightarrow{\inj_W} \Free{W} \xhookrightarrow{\iota_2}
  \Free{K +W}$  are trivially transition-preserving, since $\overline{W}$
  has no transitions and hence $\Free{\overline{W}}$ has no generating morphisms;
  \item They are coequalized by $\upsilon_{K+W}$. Indeed we have:
  \begin{align*}
    \inj_W \Cp w \Cp \iota_1 \Cp \upsilon_{K+W} &=
    \inj_W \Cp w \Cp \iota_1 \Cp [\FunK,w \Cp \FunK]\\
    &= \inj_W \Cp w \Cp \FunK\\
    &= \inj_W \Cp \iota_2 \Cp [\FunK,w \Cp \FunK]\\
    &= \inj_W \Cp \iota_2 \Cp \upsilon_{K+W}
  \end{align*}
  \item $\coeq$ is their coequalizer.
\end{itemize}
This identification of a coproduct is nothing more
than the pushout of Definition~\ref{def: addition of generators},
now seen as the result of two separate steps -- generator addition
and place merging.
\begin{equation*}
  \scalebox{0.75}{
  \begin{tikzpicture}[node distance=4cm,>=stealth',bend angle=45,auto]
    \node (1) at (0,0) {$\Free{\overline{W}}$};
    \node (2) [right = 1.25cm of 1] {$\Free{W}$};
    \node (3) [right = 1.25cm of 2] {$\Free{K}$};
    \node (4) [below = 1.5cm of 1] {$\Free{W}$};
    % \node (5) [below = 1.5cm of 2] {};
    \node (6) [below = 1.5cm of 3] {$\Free{K+W}$};
    \node (7) [below right = 1.5cm and 3cm of 6.center] {$\Free{M}$};
    \node (8) [below right = 3cm and 1.5cm of 6.center] {$\Semantics$};
    \draw[right hook-latex] (1) to node[font=\tiny] {$\inj_W$} (2);
    \draw[left hook-latex] (1) to node[font=\tiny,swap] {$\inj_W$} (4);
    \draw[->] (2) to node[font=\tiny] {$w$} (3);
    \draw[left hook-latex] (3) to node[font=\tiny] {$\iota_1$} (6);
    \draw[right hook-latex] (4) to node[font=\tiny,swap] {$\iota_2$} (6);
    \draw[transform canvas={yshift=0.5ex, xshift=0.5ex}, ->] (1) to node[font=\tiny] {$\inj_W \Cp w \Cp \iota_1$} (6);
    \draw[transform canvas={yshift=-0.5ex, xshift=-0.5ex}, ->] (1) to node[font=\tiny,swap] {$\inj_W \Cp \iota_2$} (6);
    \draw[dotted, ->] (6) to node[font=\tiny] {$\coeq$} (7);
    \draw[->] (6) to node[font=\tiny,swap] {$[\FunK, w \Cp \FunK]$} (8);
    \draw[dashed, ->] (7) to node[font=\tiny] {$\exists!$} (8);
  \end{tikzpicture}
}
\end{equation*}
\subsection{Erasing generators}
As we said, synchronizations amount to
promote computations to generators
on one hand, and to erasing other generators
on the other. We now focus on modelling what
removing generators means, starting with
some straightforward definitions:
\begin{definition}
  A net $N_w$ is called
  \emph{a subnet of the net $N$} if its places and
  transitions are a subset of places and transitions
  of $N$, and input and output functions on $N_w$
  are restrictions of the input and output functions on $N$.
\end{definition}
If $N_w$ is a subnet of $N$, then $\Free{N_w}$ is
a subcategory of $\Free{N}$, and there is an identity
on objects, identity on morphisms strict monoidal
functor, $\sub_{N_w}: \Free{N_w} \hookrightarrow \Free{N}$.
This can be used to lift the concept of subnet to $\PetriS$:
\begin{definition}
  A net $(N_w, \upsilon_{N_w})$ is a subnet of $\NetSem{N}$
  if $N_w$ is a subnet of $N$ and $\upsilon_{N_w} = \sub_{N_w} \Cp \FunN$.
\end{definition}
Notice that the subnet condition completely
specifies $\upsilon_{N_w}$ once $N_w$ is provided,
so we can also just say ``let $N_w$ be a subnet of $\NetSem{N}$''
to denote $(N_w, \upsilon_{N_w})$.
We can now use subnets and pushouts to define what the
erasing of generators is:
\begin{definition}
  Let $(N_w, \upsilon_{N_w})$ be a subnet of $\NetSem{N}$. An \emph{erasing of
  generators of $\NetSem{N}$ via $N_w$} is a subnet
  $\NetSem{K}$ of $\NetSem{N}$ such that the following
  square is a pushout, where the arrows denote the obvious inclusions:
  \begin{equation*}
    \scalebox{1}{
  \begin{tikzpicture}[node distance=2cm,>=stealth',bend angle=45,auto]
    \node (1) at (0,0) {$\Free{N}$};
    \node (2) at (4, 0) {$\Free{K}$};
    \node (3) at (0,2) {$\Free{N_w}$};
    \node (4) at (4,2) {$\Free{\overline{N_w}}$};

    \draw[left hook-latex] (4) to node[midway, above] {$\inj_{N_w}$} (3);
    \draw[->] (4) to (2);
    \draw[left hook-latex] (3) to node[midway, left] {$\sub_{N_w}$} (1);
    \draw[left hook-latex] (2) to node[midway, below] {$\sub_{K}$} (1);
  \end{tikzpicture}
}
  \end{equation*}
\end{definition}
\begin{lemma}
  \label{lem: uniqueness of erasing of generators}
  Given a net $\NetSem{N}$ and a subnet $N_w$,
  erasings of generators of $\NetSem{N}$ via
  $N_w$ are unique up to isomorphism.
\end{lemma}
Unrolling this as we did for
Definition~\ref{def: addition of generators},
$K$ represents the transitions to be added to $N_w$
to obtain $N$. $N_w$ is used to pinpoint the
generators to be erased from $N$, and $K$ represents the
``complement" of $N_w$.

We illustrate how erasing works resorting to
the example depicted in Figure~\ref{fig: Erasing generators of a net}:
In Figure~\ref{fig: example of synchronization erased net} the
generating morphisms $f$ and $g$ have been erased
from the net in Figure~\ref{fig: example of synchronization starting net 2}
via the subnet $N_w$ in Figure~\ref{fig: example of synchronization witness net}.
\begin{figure}[!ht]
  \centering
  \begin{subfigure}[t]{0.32\textwidth}\centering
    \scalebox{0.5}{
  \begin{tikzpicture}[node distance=1.3cm,>=stealth',bend angle=45,auto]
		% First net
		%
		\node [place] (1a) at (0,-1.5) {$A$};
		\node [transition] (2a) at (1.5,1) {$f$}
	  	edge [pre] (1a);
		\node [place] (3a) at (3,1) {$B$}
			edge [pre] (2a);
    \node [transition] (3b) at (3,-4)  {$h$}
     edge [pre] (1a);
    \node [transition] (4a) at (4.5,1)     {$g$}
      edge [pre] (3a);
    \node [place] (5a) at (6,1)  {$C$}
     edge [pre] (4a);
    \node [place] (5b) at (6,-4)    {$D$}
      edge [pre] (3b);
		\begin{pgfonlayer}{background}
			\filldraw [line width=4mm,join=round,\backgrnd]
				(-0.5,1.5) rectangle (6.5,-4.5);
		\end{pgfonlayer}
  \end{tikzpicture}
  }
    \caption{A net $N$.}
    \label{fig: example of synchronization starting net 2}
  \end{subfigure}
  \begin{subfigure}[t]{0.32\textwidth}\centering
    \scalebox{0.5}{
  \begin{tikzpicture}[node distance=1.3cm,>=stealth',bend angle=45,auto]
		% First net
		%
		\node [place] (1a) at (0,-1.5) {$A$};
		\node [transition] (2a) at (1.5,1) {$f$}
	  	edge [pre] (1a);
		\node [place] (3a) at (3,1) {$B$}
			edge [pre] (2a);
    \node [transition] (4a) at (4.5,1)     {$g$}
      edge [pre] (3a);
    \node [place] (5a) at (6,1)  {$C$}
     edge [pre] (4a);
		\begin{pgfonlayer}{background}
			\filldraw [line width=4mm,join=round,\backgrnd]
				(-0.5,1.5) rectangle (6.5,-4.5);
		\end{pgfonlayer}
  \end{tikzpicture}
  }
    \caption{The subnet $N_w$.}
    \label{fig: example of synchronization witness net}
  \end{subfigure}
  \begin{subfigure}[t]{0.32\textwidth}\centering
    \scalebox{0.5}{
  \begin{tikzpicture}[node distance=1.3cm,>=stealth',bend angle=45,auto]
		% First net
		%
		\node [place] (1a) at (0,-1.5) {$A$};
		\node [place] (3a) at (3,1) {$B$};
    \node [transition] (3b) at (3,-4)  {$h$}
     edge [pre] (1a);
    \node [place] (5a) at (6,1)  {$C$};
    \node [place] (5b) at (6,-4)    {$D$}
      edge [pre] (3b);
		\begin{pgfonlayer}{background}
			\filldraw [line width=4mm,join=round,\backgrnd]
				(-0.5,1.5) rectangle (6.5,-4.5);
		\end{pgfonlayer}
  \end{tikzpicture}
}
    \caption{Generators erased.}
    \label{fig: example of synchronization erased net}
  \end{subfigure}
  \caption{Erasing generators of a net.}
    \label{fig: Erasing generators of a net}
\end{figure}
\subsection{Putting things together}
Now we have the necessary concepts to characterize
all the steps happening in a synchronization. Before
proceeding, there is just a last bit of tooling we
need to develop to determine
a subnet $N_w$ from a synchronization witness $W,w$
so that the needed erasing of generators can be performed.
\begin{definition}
  Given a morphism $f$ in a FSSMC $\CategoryC$, a
  \emph{decomposition of $f$} consists in a set of morphism
  generators $\{f_1, \dots, f_n \}$ such that
  $f$ is equal to a combination of
  compositions and monoidal products of
  identities, symmetries and the $f_i$, with
  each $f_i$ used at least once. We say that each
  $f_i$ \emph{belongs} to the decomposition
  $\{f_1, \dots, f_n \}$ .
\end{definition}
\begin{example}
  Given morphism generators $f,g$ in a FSSMC, $\{ f,g \}$
  is a decomposition of $f \otimes g$, $f \otimes f \otimes g$ and $f \Cp g$.
  $f$ belongs to such decomposition. If $h$ is another morphism
  generator, then $h$ does not belong to such decomposition, and
  $\{f,g,h\}$, $\{f\}$, $\{g,h\}$ are not decompositions of $f \otimes g$ or $f \Cp g$.
\end{example}
In arbitrary categories, where further
equations between generating morphisms are
imposed, the problem of determining if
$\{ f_1, \dots, f_n \}$ constitutes a decomposition
for some morphism $f$ is connected to the word
problem, and would thus be undecidable -- think
about the fact that any group is a category, and
the word problem is undecidable for groups. Luckily,
working with FSSMCs makes the situation easier.
\begin{lemma}\label{lem: decompositions are invariant for belonging}
  Let $\CategoryC$ be a FSSMC. Every morphism admits
  exactly one decomposition.
\end{lemma}
\begin{definition}
  For a generating morphism $f$ and a morphism $g$
  in a FSSMC, we will write $f \in g$ if $f$ belongs
  to the unique decomposition of $g$.
\end{definition}
We now employ decompositions to determine a subnet
from a synchronization witness, which will determine
the generators to be erased from a net in a synchronization.
\begin{definition}
  Let $\NetSem{N}$ be a net, and let $W, w$
  be a synchronization witness for $\NetSem{N}$.
  We denote with $N_w$ the subnet of $N$ defined as follows:
  \begin{itemize}
    \item Places of $N_w$ are the generating
    objects of $\Free{N}$ that are in the image of $w$;
    \item Transitions of $N_w$ are the generating
    morphisms $g$ such that $g \in wf$ for some
    generating morphism $f$ of $\Free{W}$.
  \end{itemize}
  By definition, the functor $w:\Free{W} \to \Free{N}$ factorizes
  through $\Free{W} \xrightarrow{w'} \Free{N_w} \xhookrightarrow{\sub_{N_w}} \Free{N}$,
  where $w'$ is defined as the restriction of $w$ to $\Free{N_w}$.
\end{definition}
We conclude this section by giving the definition
of synchronization, which we worked so hard to obtain:
\begin{definition}
  \label{def: synchronization}
  Let $\NetSem{N}$ be a net, and let $W, w$
  be a synchronization witness for $\NetSem{N}$.
  A net $\NetSem{M}$ is a \emph{synchronization} of
  $\NetSem{N}$ via $W,w$ if both the following left
  and right squares are pushout squares
  for some subnet $\NetSem{K}$ of $\NetSem{N}$,
  and $\FunM$ is determined by the universal property
  of the pushout:
  \begin{equation*}
    % \scalebox{1}{
  \begin{tikzpicture}[node distance=2cm,>=stealth',bend angle=45,auto]
    \node (1) at (0,0) {$\Free{N}$};
    \node (2) at (4, 0) {$\Free{K}$};
    \node (3) at (6,0) {$\Free{M}$};
    \node (4) at (0,2) {$\Free{N_w}$};
    \node (5) at (2,2) {$\Free{W}$};
    \node (6) at (4,2) {$\Free{\overline{W}}$};
    \node (7) at (6,2) {$\Free{W}$};
    \node (8) at (8,-2) {$\Semantics$};

    \draw[left hook-latex] (4) to node[font=\tiny,midway, left] {$\sub_{N_W}$} (1);
    \draw[->] (6) to (2);
    \draw[->] (2) to (3);
    \draw[->] (7) to (3);
    \draw[->] (5) to node[font=\tiny,midway, above] {$w'$} (4);

    \draw[left hook-latex] (6) to node[font=\tiny,midway, above] {$\inj_{W}$} (5);
    \draw[right hook-latex] (6) to node[font=\tiny,midway, above] {$\inj_{W}$} (7);
    \draw[left hook-latex] (2) to node[font=\tiny,midway, below] {$\sub_{K}$} (1);

    \draw[->] (5.south) to node[font=\tiny,midway, above left] {$w$} (1.north east);

    \draw[->, bend right=10] (1) to node[font=\tiny,midway, below left] {$\FunN$} (8);
    \draw[->, bend right=10] (2) to node[font=\tiny,midway, below left] {$\FunK$} (8);
    \draw[->, bend left=10] (7) to node[font=\tiny,midway, above right] {$w \Cp \FunN$} (8);
    \draw[->, dashed] (3) to node[font=\tiny,midway, above right] {$\FunM$} (8);

  \end{tikzpicture}
% }
% https://q.uiver.app/?q=WzAsOCxbMCwwLCJcXEZyZWV7Tl93fSJdLFsxLDAsIlxcRnJlZXtXfSJdLFsyLDAsIlxcRnJlZXtcXG92ZXJsaW5lIFd9Il0sWzMsMCwiXFxGcmVle1d9Il0sWzAsMSwiXFxGcmVle059Il0sWzIsMSwiXFxGcmVle0t9Il0sWzMsMSwiXFxGcmVle019Il0sWzQsMiwiXFxTZW1hbnRpY3MiXSxbMiwzXSxbMyw2XSxbMiw1XSxbNSw2XSxbMiwxXSxbMSw0XSxbMSwwXSxbMCw0XSxbNSw0XSxbNiw3XSxbMyw3LCIiLDEseyJjdXJ2ZSI6LTF9XSxbNSw3LCIiLDEseyJjdXJ2ZSI6MX1dLFs0LDcsIiIsMSx7ImN1cnZlIjoyfV1d
% \begin{tikzcd}
% 	{\Free{N_w}} & {\Free{W}} & {\Free{\overline W}} & {\Free{W}} \\
% 	{\Free{N}} && {\Free{K}} & {\Free{M}} \\
% 	&&&& \Semantics
% 	\arrow[from=1-3, to=1-4]
% 	\arrow[from=1-4, to=2-4]
% 	\arrow[from=1-3, to=2-3]
% 	\arrow[from=2-3, to=2-4]
% 	\arrow[from=1-3, to=1-2]
% 	\arrow[from=1-2, to=2-1]
% 	\arrow[from=1-2, to=1-1]
% 	\arrow[from=1-1, to=2-1]
% 	\arrow[from=2-3, to=2-1]
% 	\arrow[from=2-4, to=3-5]
% 	\arrow[curve={height=-6pt}, from=1-4, to=3-5]
% 	\arrow[curve={height=6pt}, from=2-3, to=3-5]
% 	\arrow[curve={height=12pt}, from=2-1, to=3-5]
% \end{tikzcd}
  \end{equation*}
  In other words, $\NetSem{M}$ is a synchronization of
  $\NetSem{N}$ via $W, w$ if $\Free{M}$ is the result
  of applying the double pushout rewrite rule
  $\Free{N_w} \xleftarrow{w'} \Free{W}
  \xhookleftarrow{\inj_W} \Free{\overline{W}} \xhookrightarrow{\inj_W} \Free{W}$
  to $\Free{N}$.
\end{definition}
An example of synchronization is displayed in
Figure~\ref{fig: an example of synchronization},
where the net in
Figure~\ref{fig: example of synchronization starting net 3}
has the transitions $f, g$ synchronized to obtain
the net in
Figure~\ref{fig: example of synchronization synchronized net}.
\begin{figure}[!ht]
  \centering
  \begin{subfigure}[t]{0.49\textwidth}\centering
    \scalebox{0.5}{
  \begin{tikzpicture}[node distance=1.3cm,>=stealth',bend angle=45,auto]
		% First net
		%
		\node [place] (1a) at (0,-1.5) {$A$};
		\node [transition] (2a) at (1.5,1) {$f$}
	  	edge [pre] (1a);
		\node [place] (3a) at (3,1) {$B$}
			edge [pre] (2a);
    \node [transition] (3b) at (3,-4)  {$h$}
     edge [pre] (1a);
    \node [transition] (4a) at (4.5,1)     {$g$}
      edge [pre] (3a);
    \node [place] (5a) at (6,1)  {$C$}
     edge [pre] (4a);
    \node [place] (5b) at (6,-4)    {$D$}
      edge [pre] (3b);
		\begin{pgfonlayer}{background}
			\filldraw [line width=4mm,join=round,\backgrnd]
				(-0.5,1.5) rectangle (6.5,-4.5);
		\end{pgfonlayer}
  \end{tikzpicture}
  }
    \caption{A net $N$.}
    \label{fig: example of synchronization starting net 3}
  \end{subfigure}
  \begin{subfigure}[t]{0.49\textwidth}\centering
    \scalebox{0.5}{
  \begin{tikzpicture}[node distance=1.3cm,>=stealth',bend angle=45,auto]
		% First net
		%
		\node [place] (1a) at (0,-1.5) {$A$};
		\node [place] (3a) at (3,1) {$B$};
    \node [transition] (3b) at (3,-4)  {$h$}
     edge [pre] (1a);
    \node [transition] (3c) at (3,-1.5)     {$f \Cp g$}
      edge [pre] (1a);
    \node [place] (5a) at (6,1)  {$C$}
     edge [pre] (3c);
    \node [place] (5b) at (6,-4)    {$D$}
      edge [pre] (3b);
		\begin{pgfonlayer}{background}
			\filldraw [line width=4mm,join=round,\backgrnd]
				(-0.5,1.5) rectangle (6.5,-4.5);
		\end{pgfonlayer}
  \end{tikzpicture}
}
    \caption{Synchronizing $f$ and $g$.}
    \label{fig: example of synchronization synchronized net}
  \end{subfigure}
  \caption{An example of synchronization.}
    \label{fig: an example of synchronization}
\end{figure}
\subsection{Erasing redundant places}
As one can see, the net in
Figure~\ref{fig: example of synchronization synchronized net}
presents some isolated places.
It is obvious that these isolated places are useless if
not harmful -- e.g. in an implementation context they
constitute space-consuming information that adds up
as one performs synchronizations recursively --
and so one would like to get rid of them.
This is obtained by performing a last
step, which involves defining a new category.
\begin{definition}
  Given a FSSMC $\Free{M}$, a \emph{reduction}
  of $\Free{M}$ is a free symmetric strict monoidal
  subcategory $\Free{R} \subseteq \Free{M}$ such
  that every morphism generator of $\Free{M}$ is
  also in $\Free{M}$.

  Similarly,  $\NetSem{R}$ is is a reduction
  of $\NetSem{M}$ if $\Free{R}$ is a reduction
  of $\Free{M}$ and $\Fun{R}$ is obtained
  precomposing $\FunM$ with the obvious inclusion.
\end{definition}
\begin{lemma}
  Given a net $\NetSem{M}$, reductions of $\NetSem{M}$
  form the objects of a category, denoted $\Red{\NetSem{M}}$,
  the morphisms being the obvious inclusions. $\NetSem{M}$ is
  terminal in this category.
\end{lemma}
Notice how the conditions imposed on objects of
$\Red{\NetSem{M}}$ are very strict: Basically nothing
can be changed on morphisms, and this clearly brings
down to zero also the flexibility one has in manipulating
the objects that are in the domain or codomain
of a morphism. Hence, objects in $\Red{\NetSem{M}}$
differ between each other only with respect to how
they treat objects not connected to any morphism,
that is, how they treat objects generated by the
isolated places of $M$. But there is very little
that one can do with such objects: Indeed,
either they are included, or they aren't!

The result is that objects in $\Red{\NetSem{M}}$
only differ about how many of the isolated places of
$M$ they include as generators. Since we want to
get rid of them, we want the object of $\Red{\NetSem{M}}$
that includes none of them.
\begin{definition}
  Given a net $\NetSem{N}$, a net $\NetSem{M}$ is a
  \emph{reduced synchronization} of $\NetSem{N}$ via $W, w$ if
  $\NetSem{M}$ is initial in $\Red{\NetSem{M'}}$,
  with $\NetSem{M'}$ a synchronization of $\NetSem{N}$
  via $W, w$.
\end{definition}
Reduced synchronizations finally implement our end goal,
as depicted in
Figure~\ref{fig: five stages of a synchronization}.
To summarize, in this long explanation it has been
shown how synchronizations are not really a
single operation, but multiple operations performed in sequence.
\begin{figure}[!ht]
  \centering
  \begin{subfigure}[t]{0.32\textwidth}\centering
    \scalebox{0.5}{
  \begin{tikzpicture}[node distance=1.3cm,>=stealth',bend angle=45,auto]
		% First net
		%
		\node [place] (1a) at (0,-1.5) {$A$};
		\node [transition] (2a) at (1.5,1) {$f$}
	  	edge [pre] (1a);
		\node [place] (3a) at (3,1) {$B$}
			edge [pre] (2a);
    \node [transition] (3b) at (3,-4)  {$h$}
     edge [pre] (1a);
    \node [transition] (4a) at (4.5,1)     {$g$}
      edge [pre] (3a);
    \node [place] (5a) at (6,1)  {$C$}
     edge [pre] (4a);
    \node [place] (5b) at (6,-4)    {$D$}
      edge [pre] (3b);
		\begin{pgfonlayer}{background}
			\filldraw [line width=4mm,join=round,\backgrnd]
				(-0.5,1.5) rectangle (6.5,-4.5);
		\end{pgfonlayer}
  \end{tikzpicture}
  }
    \caption{A net $N$.}
  \end{subfigure}
  \begin{subfigure}[t]{0.32\textwidth}\centering
    \scalebox{0.5}{
  \begin{tikzpicture}[node distance=1.3cm,>=stealth',bend angle=45,auto]
		% First net
		%
		\node [place] (1a) at (0,-1.5) {$A$};
    \node [transition] (3c) at (3,-1.5)     {$f \Cp g$}
      edge [pre] (1a);
    \node [place] (5a) at (6,1)  {$C$}
     edge [pre] (3c);
		\begin{pgfonlayer}{background}
			\filldraw [line width=4mm,join=round,\backgrnd]
				(-0.5,1.5) rectangle (6.5,-4.5);
		\end{pgfonlayer}
  \end{tikzpicture}
}
    \caption{A synchronization witness $W$.}
  \end{subfigure}
  \begin{subfigure}[t]{0.32\textwidth}\centering
    \scalebox{0.5}{
  \begin{tikzpicture}[node distance=1.3cm,>=stealth',bend angle=45,auto]
		% First net
		%
		\node [place] (1a) at (0,-1.5) {$A$};
		\node [transition] (2a) at (1.5,1) {$f$}
	  	edge [pre] (1a);
		\node [place] (3a) at (3,1) {$B$}
			edge [pre] (2a);
    \node [transition] (4a) at (4.5,1)     {$g$}
      edge [pre] (3a);
    \node [place] (5a) at (6,1)  {$C$}
     edge [pre] (4a);
		\begin{pgfonlayer}{background}
			\filldraw [line width=4mm,join=round,\backgrnd]
				(-0.5,1.5) rectangle (6.5,-4.5);
		\end{pgfonlayer}
  \end{tikzpicture}
  }
    \caption{The subnet $N_w$.}
  \end{subfigure}

  \bigskip
  \begin{subfigure}[t]{0.32\textwidth}\centering
    \scalebox{0.5}{
  \begin{tikzpicture}[node distance=1.3cm,>=stealth',bend angle=45,auto]
		% First net
		%
		\node [place] (1a) at (0,-1.5) {$A$};
		\node [place] (3a) at (3,1) {$B$};
    \node [transition] (3b) at (3,-4)  {$h$}
     edge [pre] (1a);
    \node [place] (5a) at (6,1)  {$C$};
    \node [place] (5b) at (6,-4)    {$D$}
      edge [pre] (3b);
		\begin{pgfonlayer}{background}
			\filldraw [line width=4mm,join=round,\backgrnd]
				(-0.5,1.5) rectangle (6.5,-4.5);
		\end{pgfonlayer}
  \end{tikzpicture}
}
    \caption{Erasing generators ($K$).}
  \end{subfigure}
  \begin{subfigure}[t]{0.32\textwidth}\centering
    \scalebox{0.5}{
  \begin{tikzpicture}[node distance=1.3cm,>=stealth',bend angle=45,auto]
		% First net
		%
		\node [place] (1a) at (0,-1.5) {$A$};
		\node [place] (3a) at (3,1) {$B$};
    \node [transition] (3b) at (3,-4)  {$h$}
     edge [pre] (1a);
    \node [transition] (3c) at (3,-1.5)     {$f \Cp g$}
      edge [pre] (1a);
    \node [place] (5a) at (6,1)  {$C$}
     edge [pre] (3c);
    \node [place] (5b) at (6,-4)    {$D$}
      edge [pre] (3b);
		\begin{pgfonlayer}{background}
			\filldraw [line width=4mm,join=round,\backgrnd]
				(-0.5,1.5) rectangle (6.5,-4.5);
		\end{pgfonlayer}
  \end{tikzpicture}
}
    \caption{Adding generators ($M'$).}
  \end{subfigure}
  \begin{subfigure}[t]{0.32\textwidth}\centering
    \scalebox{0.5}{
  \begin{tikzpicture}[node distance=1.3cm,>=stealth',bend angle=45,auto]
		% First net
		%
		\node [place] (1a) at (0,-1.5) {$A$};
    \node [transition] (3b) at (3,-4)  {$h$}
     edge [pre] (1a);
    \node [transition] (3c) at (3,-1.5)     {$f \Cp g$}
      edge [pre] (1a);
    \node [place] (5a) at (6,1)  {$C$}
     edge [pre] (3c);
    \node [place] (5b) at (6,-4)    {$D$}
      edge [pre] (3b);
		\begin{pgfonlayer}{background}
			\filldraw [line width=4mm,join=round,\backgrnd]
				(-0.5,1.5) rectangle (6.5,-4.5);
		\end{pgfonlayer}
  \end{tikzpicture}
}
    \caption{Erasing isolated places ($M$).}
  \end{subfigure}
  \caption{Anatomy of a synchronization.}
    \label{fig: five stages of a synchronization}
\end{figure}
\subsection{Implementation perspective}
As in the case of identifications, we now look
at ways of characterizing synchronizations that make
them easier to compute algorithmically.
Addition of generators is easy to characterize:
\begin{lemma}
  \label{lem: characterization of addition of generating morphisms}
  Let $\NetSem{M}$ be an addition of generating moprhisms to
  $\NetSem{K}$ via  $W, w$. If $W$ has finite places
  and transitions, then $\Free{M}$ is isomorphic to the
  FSSMC $\CategoryC$ generated as follows:
  \begin{itemize}
    \item $\CategoryC$ and $\Free{K}$ coincide on
    object generators (and hence on objects). Concisely:
    \begin{equation*}
      \GObj{\CategoryC} = \GObj{\Free{K}}
    \end{equation*}
    \item  The generating morphisms of $\CategoryC$
    are the disjoint union of the generating morphism
    of $\Free{K}$ and $\Free{W}$, where the domain
    and codomain of generators in $W$ are substituted
    with their images through $w$. Concisely:
    \begin{equation*}
      \GMor{\CategoryC} = \GMor{\Free{K}} \sqcup
      \Suchthat{\alpha: [\Source{\alpha}/\Source{w\alpha}] \to  [\Target{\alpha}/\Target{w\alpha}]}{\alpha \in \GMor{\Free{W}}}
    \end{equation*}
  \end{itemize}
  $\FunM$ is the functor
  $\Free{M} \simeq \CategoryC \xrightarrow{F} \Semantics$,
  with $F$ defined as $\FunK$ on objects and on generating
  morphisms coming from $\Free{K}$, and as $w \Cp \FunK$
  on generating morphisms coming from $\Free{W}$.
\end{lemma}
Regarding the erasing of generators, the
characterization is again easy, and backs up our intuition:
\begin{lemma}
  \label{lem: characterization of erasing morphisms}
  Let $\NetSem{K}$ be an erasing of generators of
  $\NetSem{N}$ via $N_w$, with $N$ having
  finite places and transitions. Then $\Free{K}$
  is isomorphic to the FSSMC $\CategoryD$ generated
  as follows:
  \begin{itemize}
    \item $\CategoryD$ and $\Free{N}$ coincide on
    object generators (and hence on objects). Concisely:
    \begin{equation*}
      \GObj{\CategoryD} = \GObj{\Free{N}}
    \end{equation*}
    \item  The generating morphisms of $\CategoryD$
    are the generating morphisms of $\Free{N}$
    minus the generating morphisms of $\Free{N_w}$.
    Concisely:
    \begin{equation*}
      \GMor{\CategoryD} = \GMor{\Free{N}} /
      \Suchthat{\alpha}{\alpha \in \GMor{\Free{N_w}}}
    \end{equation*}
    Where $/$ denotes the set theoretic difference.
  \end{itemize}
  $\FunK$ is the functor
  $\Free{K} \simeq \CategoryD \xrightarrow{F} \Semantics$,
  with $F$ defined as $\FunN$ on objects and on generating
  morphisms.
\end{lemma}
Putting these lemmas together,
we can finally obtain an implementationally
viable characterization of synchronizations:
\begin{lemma}\label{lem: characterization of synchronization}
  Let $\NetSem{M}$ be a synchronization of
  $\NetSem{N}$ via $W, w$. If $W$
  has finite places and transitions, then $\Free{M}$
  is isomorphic to the FSSMC $\CategoryE$ generated
  as follows:
  \begin{itemize}
    \item $\CategoryE$ and $\Free{N}$ coincide on
    object generators (and hence on objects). Concisely:
    \begin{equation*}
      \GObj{\CategoryE} = \GObj{\Free{N}}
    \end{equation*}
    \item  The generating morphisms of $\CategoryE$
    are the generating morphisms of $\Free{N}$, to which
    are added the morphisms of $\Free{W}$ as in
    Lemma~\ref{lem: characterization of addition of generating morphisms},
    and to which are stripped away the generating
    morphisms of $\Free{N}$ that are in the decomposition of the image
    of generating morphisms of $\Free{W}$ through $w$.
    Concisely:
    \begin{align*}
      \GMor{\CategoryE} &=
      \GMor{\CategoryC} \sqcup \\
      & \sqcup \Suchthat{\alpha: [\Source{\alpha}/\Source{w\alpha}] \to  [\Target{\alpha}/\Target{w\alpha}]}{\alpha \in \GMor{\Free{W}}} / \\
      & /  \Suchthat{g \in \GMor{\Free{N}}}{\exists f \in \GMor{\Free{W}}.(g \in wf)}
    \end{align*}
    Where $/$ denotes the set theoretic difference.
  \end{itemize}
  $\FunM$ is the functor
  $\Free{M} \simeq \CategoryD \xrightarrow{F} \Semantics$,
  with $F$ defined as $\FunN$ on objects and on generating
  morphisms coming from $\Free{N}$, and as $w \Cp \FunN$
  on generating morphisms coming from $\Free{W}$.
\end{lemma}
As for identifications, these lemmas have
the advantage of performing synchronizations
essentially by means of substitutions, which are
computationally more efficient than calculating
pushouts using quotients.
\section{Recovering old concepts: Gluing between nets}
  \label{sec: recovering old concepts gluing between nets}
In this section, we show how the formalism
of ``gluing places/transitions within a net''
can specialize to the usual ``gluings between different nets''
formalism already covered in the literature.
We start by noticing that
the monoidal structure of $\PetriS$ describes the
operation of ``putting nets next to each other''.
With this, we can further specialize
the concepts defined in Section~\ref{sec: identifications}.
\begin{remark}
  Consider a couple of Petri nets
  $\NetSem{M}, \NetSem{N}$. Taking a net $W$ and
  transition-preserving functors
  \begin{equation*}
    l: \Free{W} \to \Free{M} \qquad r: \Free{W} \to \Free{N}
  \end{equation*}
  Such that  $l \Cp \FunM = r \Cp \FunN$, combining $l, r$ with the
  canonical injections of the coproduct one obtains
  \begin{equation*}
    l': \Free{W} \to \Free{M+N} \qquad r': \Free{W} \to \Free{N+N}
  \end{equation*}
  From this, using the properties of the
  coproduct one infers that $[\FunM, \FunN]$
  coequalizes $l', r'$.
  Applying Lemma~\ref{lem: characterizing coequalizers}
  the coequalizer of $l', r'$ can be calculated,
  obtaining an identification for the net $\NetSem{(N + M)}$.
  Since the coequalizer of the coproduct
  is nothing more than the pushout,
  we have recovered the familiar gluing
  of nets developed in~\cite{Fong2016} using the formalism
  of~\cite{Fong2015}.
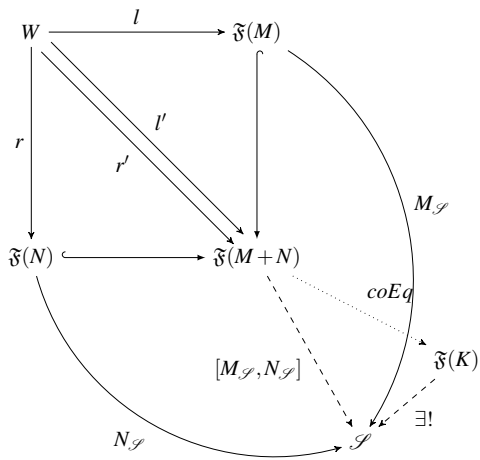
  \begin{equation*}
    \scalebox{0.75}{
    \begin{tikzpicture}[node distance=4cm,>=stealth',bend angle=45,auto]
        \node (1) at (0,0) {$W$};
        \node (2) [right of=1] {$\Free{M}$};
        \node (3) [below of=1] {$\Free{N}$};
        \node (4) [below of=2] {$\Free{M+N}$};
        \node (5) [below right = 3cm and 1.5cm of 4.center] {$\Semantics$};
        \node (6) [below right = 1.5cm and 3cm of 4.center] {$\Free{K}$};

        \draw[->] (1) to node {$l$} (2);
        \draw[->] (1) to node[swap] {$r$} (3);

        \draw[right hook-latex] (2) to node {} (4);
        \draw[right hook-latex] (3) to node {} (4);

        \draw[transform canvas={yshift=0.5ex, xshift=0.5ex}, ->] (1) to node {$l'$} (4);
        \draw[transform canvas={yshift=-0.5ex, xshift=-0.5ex}, ->] (1) to node[swap] {$r'$} (4);

        \draw[bend left, ->] (2) to node {$\FunM$} (5);
        \draw[bend right, ->] (3) to node[swap] {$\FunN$} (5);

        \draw[dashed, ->] (4) to node[swap] {$[\FunM, \FunN]$} (5);
        \draw[dotted, ->] (4) to node {$\coeq$} (6);
        \draw[dashed, ->] (6) to node {$\exists !$} (5);
    \end{tikzpicture}
}
  \end{equation*}
\end{remark}
We conclude this section by giving an intuitive explanation of how
net composition as in~\cite{Bruni2013} can be recaptured as
a special case of synchronization in our framework. Pinning down
the right definitions to isolate these synchronizations as a subclass of
all the available ones has been proven lenghty and elusive, and would
be beyond the scope of this paper. We are confident that the reader will be
satisfied with this conceptual explanation, and will be able to fill in
the needed details as a (long and unrewarding) exercise.
\begin{remark}
  In~\cite{Bruni2013}, nets come endowed with left and right ports,
  as depicted in Figure~\ref{fig: gluing nets with boundaries pawel}.
  \begin{figure}[!ht]
    \centering
    \begin{subfigure}[t]{0.49\textwidth}\centering
     \scalebox{0.5}{
  \begin{tikzpicture}[node distance=1.2cm,>=stealth',bend angle=45,auto]
		% First net
		%
		\node [place,tokens=0] (1a) at (0,0) {};
    \node [transition] (2a) at (1.5,0)     {}
    edge [pre] node[swap] {$2$} (1a);
		\node [minimum size = 2pt, circle, fill = black] (3a) at (3,0)  {}
          edge [pre] (2a);
		\node [place,tokens=0] (3b) at (1, -1) {}
          edge [pre] (2a);
    \node (01a) [right of = 3a] {};
    \node [minimum size = 2pt, circle, fill = black] (11a) at (4,0){};
    \node [transition] (21a) at (5.5,1)    {}
    edge [pre] node[swap] {$2$} (11a);
    \node [transition] (21b) at (5.5, -1)     {}
    edge [pre] (11a);
    \node [place,tokens=0] (31a) at (6.5, 1) {}
          edge [pre] (21a);
    \node [place,tokens=0] (31b) at (6.5, -1) {}
    edge [pre] (21b);
    \begin{pgfonlayer}{background}
			\filldraw [line width=10mm,join=round,\backgrnd]
        (-0.25,3) rectangle (6.75,-3);
      \draw [line width=1pt,join=round,black]
        (-0.5,1.5) rectangle (3,-1.5);
      \draw [line width=1pt,join=round,black]
        (4,1.5) rectangle (7,-1.5);
		\end{pgfonlayer}
  \end{tikzpicture}
}
      \caption{Two nets with boundaries.}
      \label{fig: pawel basic net}
    \end{subfigure}
    \begin{subfigure}[t]{0.49\textwidth}\centering
      \scalebox{0.5}{
  \begin{tikzpicture}[node distance=1.2cm,>=stealth',bend angle=45,auto]
		% First net
		%
		\node [place,tokens=0] (1a) at (0,0) {};
    \node [transition] (2a) at (3.5,1)  {}
    edge [pre] node[swap]{4}(1a);

    \node [transition] (2a') at (3.5,-1)  {}
    edge [pre] node{2}(1a);
		\node [place,tokens=0] (3b) at (2.5,0) {}
          edge [pre] node[swap] {$2$} (2a)
          edge [pre] (2a');
    \node [place,tokens=0] (31a) at (6.5,1) {}
          edge [pre] (2a);
    \node [place,tokens=0] (31b) at (6.5,-1)  {}
    edge [pre] (2a');   
    \begin{pgfonlayer}{background}
			\filldraw [line width=10mm,join=round,\backgrnd]
        (-0.25,3) rectangle (6.75,-3);
      \draw [line width=1pt,join=round,black]
        (-0.5,1.5) rectangle (7,-1.5);
		\end{pgfonlayer}
  \end{tikzpicture}
  }
      \caption{Their identification.}
      \label{fig: pawel nets fused}
    \end{subfigure}
      \caption{Gluing nets with boundaries as in~\cite{Bruni2013}.}
      \label{fig: gluing nets with boundaries pawel}
  \end{figure}

  \noindent
  As these nets do not come endowed with
  any semantics, we did not decorate places and transitions.
  They are sorted by the number of left and right
  ports they posses: A net having $m$ left ports and $n$
  right ports is interpreted as a morphsim $m \to n$.
  Nets and ports arranged this way form a
  monoidal category, and nets are composed by connecting
  them along their ports. When nets are composed,
  a suitable (minimal) multiset of transitions from the
  net on the left is synchronized with a corresponding suitable (minimal)
  multiset of transitions of the net on the right.
  Inputs and outputs for the new transitions are calculated 
  by requiring that all the tokens produced by the transitions of the
  net on the left are consumed by the transitions of the net on the right.

  This is better understood by looking 
  at Figure~\ref{fig: pawel basic net},  
  where there are two transitions connected to the left port 
  in the right net; the top one is expecting two tokens from
  the port, whereas the bottom one is expecting one.
  On the other hand the left net consists of only one transition
  connected to its right port, producing 
  one token per firing.
  These nets are composable, and the result is shown in Figure~\ref{fig: pawel nets fused}:
  The transitions in the resulting composition represent a conflation of events:
  \begin{itemize}
    \item The first transition represents the event of the 
    transition in the left net firing twice, followed by 
    the firing of the top transition of the right net;
    \item The second transition represents the event of the transition
    in the left net firing once, followed by the firing of the bottom
    transition in the right net.
  \end{itemize}
  Evidently, domains and codomains of these transitions are adjusted accordingly:
  For instance, the topmost transition consumes four tokens, corresponding to
  the tokens consumed by the original transition in the left net multiplied by the number
  of time it has to fire.

  Interpreting this kind of synchronization 
  in our framework is complicated by the fact that to-be-syncronized transitions
  in the left and right net must be suitably copied to satisfy the
  minimal multiset condition. Moreover, we need to find a 
  suitable substitute for ports, which are not 
  available to us. So, as a first step, we replace ports with 
  places. This is consistent with interpreting 
  Petri net transitions as processes, which 
  consume and produce resources of a given type.
  The net of Figure~\ref{fig: pawel basic net} is
  shown with places replacing ports in
  Figure~\ref{fig: pawel nets in our framework}.
  \begin{figure}[!ht]
    \centering
    \begin{subfigure}[t]{0.34\textwidth}\centering
      \scalebox{0.5}{
\begin{tikzpicture}[node distance=1.2cm,>=stealth',bend angle=45,auto]
  % First net
  %
  \node [place,tokens=0] (1a) at (0,0) {$A$};
  \node [transition] (2a) at (1.5,0)   {$f$}
  edge [pre] node[swap] {$2$} (1a);
  \node [place,tokens=0] (3a) at (3,0)  {$C$}
        edge [pre] (2a);
  \node [place,tokens=0] (3b) at (1,-1)  {$B$}
        edge [pre] (2a);
  \node [place,tokens=0] (11a) at (4,0)  {$C$};
  \node [transition] (21a) at (5.5,1) {$h$}
  edge [pre] node[swap] {$2$} (11a);
  \node [transition] (21b) at (5.5,-1)  {$k$}
  edge [pre] (11a);
  \node [place,tokens=0] (31a) at (6.5,1) {$D$}
        edge [pre] (21a);
  \node [place,tokens=0] (31b) at (6.5,-1) {$E$}
  edge [pre] (21b);
  \begin{pgfonlayer}{background}
    \filldraw [line width=10mm,join=round,\backgrnd]
      (-0.25,3) rectangle (6.75,-3);
    \draw [line width=1pt,join=round,black]
      (-0.5,1.5) rectangle (3,-1.5);
    \draw [line width=1pt,join=round,black]
      (4,1.5) rectangle (7,-1.5);
  \end{pgfonlayer}
\end{tikzpicture}
}
      \caption{The same nets, with semantics.}
      \label{fig: pawel nets in our framework}
    \end{subfigure}
    \begin{subfigure}[t]{0.34\textwidth}\centering
      \scalebox{0.5}{
\begin{tikzpicture}[node distance=1.2cm,>=stealth',bend angle=45,auto]
  % First net
  %
  \node [place,tokens=0] (1a) at (0,.5) {$A$};
  \node [transition] (2a) at (1.5,2)   {$f$}
  edge [pre] node[swap] {$2$} (1a);
  \node [place,tokens=0] (3a) at (3,0.5)  {$C$}
        edge [pre] (2a);
  \node [place,tokens=0] (3b) at (1.5,.5)  {$B$}
        edge [pre] (2a);
  %

  % Second net
  \node [place,tokens=0] (1a') at (0,-.5) {$A$};
  \node [transition] (2a') at (1.5,-2)   {$f$}
  edge [pre] node[swap] {$2$} (1a');
  \node [place,tokens=0] (3a') at (3,-0.5)  {$C$}
        edge [pre] (2a');
  \node [place,tokens=0] (3b') at (1.5,-.5)  {$B$}
        edge [pre] (2a');
  %

  % Third net
  \node [place,tokens=0] (11a) at (4,0)  {$C$};
  \node [transition] (21a) at (5.5,1) {$h$}
  edge [pre] node[swap] {$2$} (11a);
  \node [transition] (21b) at (5.5,-1)  {$k$}
  edge [pre] (11a);
  \node [place,tokens=0] (31a) at (6.5,1) {$D$}
        edge [pre] (21a);
  \node [place,tokens=0] (31b) at (6.5,-1) {$E$}
  edge [pre] (21b);   
  \begin{pgfonlayer}{background}
    \filldraw [line width=10mm,join=round,\backgrnd]
      (-0.25,3) rectangle (6.75,-3);
    \draw [line width=1pt,join=round,black]
      (-0.5,.5) rectangle (3,2.5);
    \draw [line width=1pt,join=round,black]
      (-0.5,-.5) rectangle (3,-2.5);
    \draw [line width=1pt,join=round,black]
      (4,1.5) rectangle (7,-1.5);
  \end{pgfonlayer}
\end{tikzpicture}
}
      \caption{Domain net suitably duplicated.}
      \label{fig: pawel nets duplicated}
    \end{subfigure}
    \begin{subfigure}[t]{0.34\textwidth}\centering
      \scalebox{0.5}{
\begin{tikzpicture}[node distance=1.2cm,>=stealth',bend angle=45,auto]
  % First net
  %
  \node [place,tokens=0] (1a) at (0,0) {$A$};
  \node [transition] (2a) at (1.5,-1)   {$f$}
        edge [pre] node {$2$} (1a);
  \node [transition] (2a') at (1.5,1)   {$f$}
        edge [pre] node[swap] {$2$} (1a);
  \node [place,tokens=0] (3a) at (3,0)  {$C$}
        edge [pre] (2a)
        edge [pre] (2a');
  \node [place,tokens=0] (3b) at (1.5,0)  {$B$}
        edge [pre] (2a)
        edge [pre] (2a');
  \node [place,tokens=0] (11a) at (4,0)  {$C$};
  \node [transition] (21a) at (5.5,1) {$h$}
        edge [pre] node[swap] {$2$} (11a);
  \node [transition] (21b) at (5.5,-1)  {$k$}
        edge [pre] (11a);
  \node [place,tokens=0] (31a) at (6.5,1) {$D$}
        edge [pre] (21a);
  \node [place,tokens=0] (31b) at (6.5,-1) {$E$}
        edge [pre] (21b);   
  \begin{pgfonlayer}{background}
    \filldraw [line width=10mm,join=round,\backgrnd]
      (-0.25,3) rectangle (6.75,-3);
    \draw [line width=1pt,join=round,black]
      (-0.5,1.5) rectangle (3,-1.5);
      \draw [line width=1pt] (-0.5,0) -- (3,0);
    \draw [line width=1pt,join=round,black]
      (4,1.5) rectangle (7,-1.5);
  \end{pgfonlayer}
\end{tikzpicture}
}
      \caption{Indentification of copies.}
      \label{fig: pawel nets first identification}
    \end{subfigure}
    \begin{subfigure}[t]{0.31\textwidth}\centering
      \scalebox{0.5}{
\begin{tikzpicture}[node distance=1.2cm,>=stealth',bend angle=45,auto]
  % First net
  %
  \node [place,tokens=0] (1a) at (0,0) {$A$};
  \node [transition] (2a) at (1.5,-1)   {$f$}
        edge [pre] node {$2$} (1a);
  \node [transition] (2a') at (1.5,1)   {$f$}
        edge [pre] node[swap] {$2$} (1a);
  \node [place,tokens=0] (3a) at (3.5,0)  {$C$}
        edge [pre] (2a)
        edge [pre] (2a');
  \node [place,tokens=0] (3b) at (1.5,0)  {$B$}
        edge [pre] (2a)
        edge [pre] (2a');
  \node [transition] (21a) at (5.5,1) {$h$}
  edge [pre] node[swap] {$2$} (3a);
  \node [transition] (21b) at (5.5,-1)  {$k$}
  edge [pre] (3a);
  \node [place,tokens=0] (31a) at (6.5,1) {$D$}
        edge [pre] (21a);
  \node [place,tokens=0] (31b) at (6.5,-1) {$E$}
  edge [pre] (21b);
  \begin{pgfonlayer}{background}
    \filldraw [line width=10mm,join=round,\backgrnd]
      (-0.25,3) rectangle (6.75,-3);
    \draw [line width=1pt,join=round,black]
      (-0.5,1.5) rectangle (3.5,-1.5);
    \draw [line width=1pt] (-0.5,0) -- (3,0);
    \draw [line width=1pt,join=round,black]
      (3.5,1.5) rectangle (7,-1.5);
  \end{pgfonlayer}
\end{tikzpicture}
}
      \caption{Identification of ports.}
      \label{fig: pawel nets second identification}
    \end{subfigure}
    \begin{subfigure}[t]{0.31\textwidth}\centering
       \scalebox{0.5}{
  \begin{tikzpicture}[node distance=1.2cm,>=stealth',bend angle=45,auto]
		% First net
		%
		\node [place,tokens=0] (1a) at (0,0) {$A$};
    \node [transition] (2a) at (3.5,1)  {$(f \otimes f) \Cp h$}
    edge [pre] node[swap]{4}(1a);

    \node [transition] (2a') at (3.5,-1)  {$f \Cp k$}
    edge [pre] node{2}(1a);
		\node [place,tokens=0] (3b) at (2.5,0) {$B$}
          edge [pre] node[swap] {$2$} (2a)
          edge [pre] (2a');
    \node [place,tokens=0] (31a) at (6.5,1) {$D$}
          edge [pre] (2a);
    \node [place,tokens=0] (31b) at (6.5,-1)  {$E$}
    edge [pre] (2a');   
    \begin{pgfonlayer}{background}
			\filldraw [line width=10mm,join=round,\backgrnd]
        (-0.25,3) rectangle (6.75,-3);
      \draw [line width=1pt,join=round,black]
        (-0.5,1.5) rectangle (7,-1.5);
		\end{pgfonlayer}
  \end{tikzpicture}
}
      \caption{Synchronization.}
      \label{fig: pawel nets synchronized}
    \end{subfigure}
      \caption{Gluing nets with boundaries in our formalism.}
      \label{fig: gluing nets with boundaries}
  \end{figure}

  \noindent
  Net composition is performed in four stages: First, 
  we suitably duplicate the component nets to meet the minimal
  multiset condition, as in Figure~\ref{fig: pawel nets duplicated}.
  Then we identify everything but the transitions connected to
  the to-be-composed ports, obtaining components with a suitable
  number of copies of the transitions, as in
  Figure~\ref{fig: pawel nets first identification}.
  Having done this, we are ready to identify 
  places of the two resulting component nets, as shown
  in Figure~\ref{fig: pawel nets second identification}.
  In the resulting net, to-be-merged transitions
  now share resources, and can thus be synchronized.
  The synchronization is performed exactly as
  prescribed in~\cite{Bruni2013}, finding the smallest total
  number of firings allowing no remaining tokens
  between the connections. The result is shown 
  in Figure~\ref{fig: pawel nets synchronized}, 
  where it is evident by looking at the 
  decoration that $f$ fires respectively one and 
  two times, while $h, k$ fire time each.
\end{remark}
\section{Changing semantics}
  \label{sec: changing semantics}
We now focus on what happens when
the semantics -- the category $\Semantics$ --
changes. By definition, we can see $\PetriS$
as the slice category of nets and morphisms between
their corresponding FSSMCs over $\Semantics$.
By using very well known results, it is clear then
that any functor $\Semantics_1 \to \Semantics_2$
induces a corresponding monoidal functor
$\PetriG{1} \to \PetriG{2}$. Hence we have:
\begin{lemma}
  There is a category having $\PetriS$, for each
  $\Semantics$, as objects, and strict monoidal functors
  $\Semantics_1 \to \Semantics_2$ as morphisms between
  them.
\end{lemma}
The next lemma, which follows directly from the
definitions, shows how synchronizations, identifications
and transition-preserving functors are preserved
when changing semantics.
\begin{lemma}\label{lem: semantic change morphism}
  Let $F: \Free{M} \to \Free{N}$
  be a synchronization (resp.
  transition-preserving functor, identification)
  between $(M, M_{\Semantics_1})$ and
  $(N, N_{\Semantics_1})$ in $\PetriG{1}$.
  If there is a morphism $\PetriG{1} \to \PetriG{2}$
  then $F$ is a synchronization (resp.
  transition-preserving functor, identification) between
  $(M, M_{\Semantics_2})$ and $(N, N_{\Semantics_2})$
  in $\PetriG{2}$.
\end{lemma}
Even more importantly it is worth noticing that, given
$(M, M_{\Semantics_1})$ in $\PetriG{1}$ and
$(M, M_{\Semantics_2})$ in $\PetriG{2}$, by applying
the universal property of the product we can consider
$(M, M_{\Semantics_1 \times \Semantics_2})$ in
$\Petri^{\Semantics_1 \times \Semantics_2}$,
where $M_{\Semantics_1 \times \Semantics_2}$ is
defined to be just $\langle M_{\Semantics_1},
M_{\Semantics_2} \rangle$.
If a functor $F: \Free{M} \to \Free{N}$ commutes
both with $M_{\Semantics_1}, N_{\Semantics_1}$
and with $M_{\Semantics_2}, N_{\Semantics_2}$,
then again by the universal property of the product
it commutes with $M_{\Semantics_1 \times \Semantics_2}$
and $N_{\Semantics_1 \times \Semantics_2}$,
which makes the following lemma obvious to prove:
\begin{lemma}\label{lem: semantic change product}
  Let $F: \Free{M} \to \Free{N}$
  be a synchronization (resp.
  ransition-preserving functor, identification)
  between $(M, M_{\Semantics_1})$ and
  $(N, N_{\Semantics_1})$ in $\PetriG{1}$,
  and between $(M, M_{\Semantics_2})$
  and $(N, N_{\Semantics_2})$ in $\PetriG{2}$.
  Then $F$ is a synchronization (resp.
  transition-preserving functor, identification) between
  $(M, M_{\Semantics_1 \times \Semantics_2})$ and
  $(N, N_{\Semantics_1 \times \Semantics_2})$
  in $\Petri^{\Semantics_1 \times \Semantics_2}$.
\end{lemma}
Lemmas~\ref{lem: semantic change morphism}
and~\ref{lem: semantic change product} are of
great practical value, because they make possible
the design of a semantics in a compartimentalized
fashion. For instance, having modelled the
skeleton of some problem using Petri nets~\cite{Genovese2018a},
we could have semantics
$\Semantics_1, \dots, \Semantics_n$ capturing
different computational aspects of the problem at
hand -- for instance, one such semantics may
represent guards~\cite{Genovese2020}.
Such semantics can be conbined using functors
$\Semantics_1 \times \dots \times \Semantics_n \to
\Semantics$, which are interpreted as procedures to
use compartimentalized information to
compute more complicated tasks.
Lemmas~\ref{lem: semantic change morphism}
and~\ref{lem: semantic change product} guarantee
that, in such setting, synchronizations, identifications
and transition-preserving functors that are
compatible with $\Semantics_1 \times \dots
\times \Semantics_n$ are also compatible with
$\Semantics$.
\section{Conclusion and future work}
  \label{sec: conclusion and future work}
In this work, we considered Petri nets 
endowed with a semantics as the object of our 
study, and generalized the notion of morphism 
between Petri nets allowing to map transitions 
to entire sequences of computations. As a result, 
we obtained categories $\PetriS$ which are 
parametrized over $\Semantics$. In each of such 
categories, we identified two particular classes 
of morphisms, which we called \emph{synchronizations} 
and \emph{identifications}. These can be
respectively used to model the notion of gluing nets 
along transitions and places, two techniques that 
are by now considered standard in the literature.

In our model, it becomes apparent how the main 
difference between these two well-known styles of 
gluing is that whereas merging on places acts 
at a purely topological level -- places are mapped 
to places and transitions to transitions -- 
merging on transitions is better understood as 
acting at a computational level -- places are 
mapped to places, but transitions are mapped 
to computations. On a deeper technical level,
we established how merging on places arises
as the result of a single pushout operation,
while merging on transitions is the result of a
double pushout operation.

Moreover, we emphasized how the notion of 
net gluing -- be it a synchronization or an 
identification -- is completely independent from 
the concept of ports, which has been considered 
standard in the last decade of research. 
Indeed, net gluings can be perfomed within 
the same net, and are even better understood 
this way. The familiar notion of gluing 
different nets together can be recovered by 
using the obvious coproduct structure present in 
$\PetriS$ to define a monoidal product.

Directions of future work will be mainly developed
along three different axes. The first involves 
implementing the theory hereby formalized in a 
formally verified setting, the importance of this 
goal being already highlighted by the constant attention we put 
in investigating the computational feasibility of 
the constructions we defined. 

The second involves further generalizing the 
notion of morphism in $\PetriS$. In particular, 
we believe that switching to a profunctorial setting 
may reveal some other interesting classes of morphsims 
which may have practical applications. 

The third direction of research involves finding
a suitable graphical calculus for synchronization 
and identification, which does not depend on the 
concept of ports. Regarding this, it is our opinion 
that a graphical formalism based on \emph{optics}~\cite{Pickering2017} 
may be a good point to start from.
\section*{Acknowledgements}
The author wants to thank his fellow team members at Statebox and  David Spivak for the useful discussions. David Spivak and Christina Vasilakopoulou have also to be credited for pointing out how, maybe, synchronizations could admit a description in terms of double pushout rewriting, a fact that prompted the author to rework some definitions and extensively rewrite the paper.
\bibliographystyle{plain}
\bibliography{Bibliography}

\begin{thebibliography}{10}

\bibitem{Baez}
J.~C. Baez.
\newblock Network {{Theory}} ({{Part}} 1).
\newblock available at
  \url{https://johncarlosbaez.wordpress.com/2011/03/04/network-theory-part-1/}.

\bibitem{Baez2013}
J.~C. Baez.
\newblock Quantum techniques for reaction networks.
\newblock {\em Advances in Mathematical Physics}, 2018:1--9, November 2018.

\bibitem{9470566}
J.~C. Baez, F.~Genovese, J.~Master, and M.~Shulman.
\newblock Categories of nets.
\newblock In {\em 2021 36th Annual ACM/IEEE Symposium on Logic in Computer
  Science (LICS)}, pages 1--13, Los Alamitos, CA, USA, jul 2021. IEEE Computer
  Society.

\bibitem{Baez2018}
J.~C. Baez and J.~Master.
\newblock Open {Petri} nets.
\newblock {\em Mathematical Structures in Computer Science}, 30(3):314–341,
  2020.

\bibitem{Baldan2009}
P.~Baldan, F.~Bonchi, and F.~Gadducci.
\newblock Encoding asynchronous interactions using open petri nets.
\newblock In {\em {CONCUR} 2009 - Concurrency Theory}, pages 99--114. Springer
  Berlin Heidelberg, 2009.

\bibitem{Baldan2015}
P.~Baldan, F.~Bonchi, F.~Gadducci, and Giacoma~Valentina Monreale.
\newblock Modular encoding of synchronous and asynchronous interactions using
  open petri nets.
\newblock {\em Science of Computer Programming}, 109:96--124, October 2015.

\bibitem{Baldan2003}
P.~Baldan, R.~Bruni, and U.~Montanari.
\newblock Pre-nets, {{Read Arcs}} and {{Unfolding}}: {{A Functorial
  Presentation}}.
\newblock In Martin Wirsing, Dirk Pattinson, and Rolf Hennicker, editors, {\em
  Recent {{Trends}} in {{Algebraic Development Techniques}}}, volume 2755,
  pages 145--164. {Springer Berlin Heidelberg}.

\bibitem{Bruni2013}
R.~Bruni, H.~Melgratti, U.~Montanari, and P.~Soboci\'nski.
\newblock Connector algebras for {C/E} and {P/T} nets' interactions.
\newblock {\em Logical Methods in Computer Science}, 9(3), September 2013.

\bibitem{mustcite01}
R.~Bruni, Jos{\'{e}} Meseguer, U.~Montanari, and V.~Sassone.
\newblock Functorial models for {P}etri {N}ets.
\newblock {\em Inf. Comput.}, 170(2):207--236, 2001.

\bibitem{Fong2016}
B.~Fong.
\newblock {\em The {{Algebra}} of {{Open}} and {{Interconnected Systems}}}.
\newblock Phd thesis.

\bibitem{Fong2015}
B.~Fong.
\newblock {D}ecorated {C}ospans.
\newblock {\em Theory and Applications of Categories}, 30(33):1096--1120, 2015.

\bibitem{Genovese2018a}
F.~Genovese.
\newblock Behavioral {{Programming}} with {{Petri Nets}} a la {{Functional
  Way}}: {$\mkern1mu$}{{Smart Contracts}}.

\bibitem{Genovese2019}
F.~Genovese, A.~Gryzlov, Jelle Herold, M.~Perone, E.~Post, and A.~Videla.
\newblock Computational {{Petri Nets}}: {{Adjunctions Considered Harmful}}.

\bibitem{Genovese2021hierar}
F.~Genovese, F.~Loregian, and D.~Palombi.
\newblock A categorical semantics for hierarchical petri nets.
\newblock {\em Electronic Proceedings in Theoretical Computer Science},
  350:51--68, December 2021.

\bibitem{Genovese2021mana}
F.~Genovese, F.~Loregian, and D.~Palombi.
\newblock Nets with mana: a framework for chemical reaction modelling.
\newblock In {\em Graph Transformation}, pages 185--202. Springer International
  Publishing, 2021.

\bibitem{Genovese2022bdd}
F.~Genovese, F.~Loregian, and D.~Palombi.
\newblock A categorical semantics for bounded petri nets.
\newblock {\em Electronic Proceedings in Theoretical Computer Science},
  372:59--71, November 2022.

\bibitem{Genovese2020}
F.~Genovese and David~I. Spivak.
\newblock A categorical semantics for {G}uarded {P}etri {N}ets.
\newblock pages 57--74, 2020.

\bibitem{Genovese2020guarde}
Fabrizio Genovese and David~I. Spivak.
\newblock A categorical semantics for guarded petri nets.
\newblock In {\em Graph Transformation}, pages 57--74. Springer International
  Publishing, 2020.

\bibitem{habel_muller_plump_2001}
A.~Habel, J.~Müller, and D.~Plump.
\newblock Double-pushout graph transformation revisited.
\newblock {\em Mathematical Structures in Computer Science}, 11(5):637–688,
  2001.

\bibitem{mustcite02}
P.~Katis, N.~Sabadini, and R.~F.~C. Walters.
\newblock Representing place/transition nets in {Span(Graph)}.
\newblock In Michael Johnson, editor, {\em Algebraic Methodology and Software
  Technology, 6th International Conference, {AMAST} '97, Sydney, Australia,
  December 13-17, 1997, Proceedings}, volume 1349 of {\em Lecture Notes in
  Computer Science}, pages 322--336. Springer, 1997.

\bibitem{Licata2011}
D.~Licata.
\newblock Running {{Circles Around}} ({{In}}) {{Your Proof Assistant}}; or,
  {{Quotients}} that {{Compute}}.

\bibitem{Master2019}
J.~Master.
\newblock Petri nets based on {L}awvere theories.
\newblock {\em Mathematical Structures in Computer Science}, 30(7):833–864,
  2020.

\bibitem{Meseguer1990}
J.~Meseguer and U.~Montanari.
\newblock Petri nets are monoids.
\newblock 88(2):105--155.

\bibitem{Pickering2017}
M.~Pickering, J.~Gibbons, and N.~Wu.
\newblock Profunctor optics: Modular data accessors.
\newblock {\em The Art, Science, and Engineering of Programming}, 1(2), April
  2017.

\bibitem{Sassone1995}
V.~Sassone.
\newblock On the {{Category}} of {{Petri Net Computations}}.
\newblock In Peter~D. Mosses, Mogens Nielsen, and Michael~I. Schwartzbach,
  editors, {\em {{TAPSOFT}} '95: {{Theory}} and {{Practice}} of {{Software
  Development}}}, volume 915, pages 334--348. {Springer Berlin Heidelberg}.

\bibitem{Sobocinski2013a}
P.~Soboci\'nski and Owen Stephens.
\newblock Penrose: {{Putting Compositionality}} to {{Work}} for {{Petri Net
  Reachability}}.
\newblock In Reiko Heckel and Stefan Milius, editors, {\em Algebra and
  {{Coalgebra}} in {{Computer Science}}}, volume 8089, pages 346--352.
  {Springer Berlin Heidelberg}.

\bibitem{Sobocinski2013}
P.~Soboci\'nski and Owen Stephens.
\newblock Reachability via {{Compositionality}} in {{Petri}} nets.

\bibitem{StateboxTeam2019a}
{Statebox Team}.
\newblock The {{Mathematical Specification}} of the {{Statebox Language}}.

\end{thebibliography}
\clearpage
\appendix
\section{Proofs}
\begingroup
\def\thetheoremUnified{\ref{lem: free distributes over coproducts}}
\begin{lemma}
  Given nets $M, N$, it is $\Free{M + N} \simeq \Free{M} + \Free{N}$.
\end{lemma}
\addtocounter{theoremUnified}{-1}
\endgroup
\begin{proof}
  The proof is obvious by noting that $+$ works
  by taking disjont unions both of the places/transitions
  of nets and of the generating objects/morphisms of FSSMCs.
  Hence, there is a bijection between the generating objects
  and morphisms of $\Free{M + N}$ and $\Free{M} + \Free{N}$,
  from which an isomorphism of categories can be built
  by using freeness.
\end{proof}
\begingroup
\def\thetheoremUnified{\ref{lem: monoidal structure of PetriS}}
\begin{lemma}
  $\PetriS$ inherits a monoidal structure from the
  biproduct structure in the category of symmetric monoidal categories.\end{lemma}
\addtocounter{theoremUnified}{-1}
\endgroup
\begin{proof}
    First of all we need to prove that if
    $\CategoryC$ and $\CategoryD$ are free symmetric
    strict monoidal categories, so is
    $\CategoryC + \CategoryD$. This is
    easy to prove by noticing that by definition
    that generating objects and morphisms in
    $\CategoryC + \CategoryD$ are given by taking
    the disjoing union of generating objects
    and morphisms, respectively, of $\CategoryC$ and $\CategoryD$.
    Moreover, denoting with $M+N$ the net
    obtained by taking the disjoint union, respectively,
    of places and transitions of $M, N$, we
    have already proved that
    $\Free{M + N} \simeq \Free{M} + \Free{N}$.
    From this, we can set
    $\NetSem{M} + \NetSem{N} = \NetSem{(M + N)}$,
    where $\upsilon_{M + N}$ is
    defined by applying the universal property of coproducts:
    \begin{equation*}
        % \scalebox{0.75}{
    \begin{tikzpicture}[node distance=4cm,>=stealth',bend angle=45,auto]
        \node (1) at (0,0) {$\Free{M} + \Free{N} = \Free{M+N}$};
        \node (2) [right = 2cm of 1] {$\Free{M}$};
        \node (3) [below = 1cm of 1] {$\Free{N}$};
        \node (4) [below = 1cm of 2] {$S$};

        \draw[->] (2) to node[font=\tiny] {$\FunM$} (4);
        \draw[->] (3) to node[font=\tiny,swap] {$\FunN$} (4);
        \draw[left hook-latex] (2) to node[font=\tiny] {} (1);
        \draw[right hook-latex] (3) to node[font=\tiny] {} (1);
        \draw[dashed, ->] (1) to node[font=\tiny] {$[\FunM, \FunN]$} (4);
      \end{tikzpicture}
% }
    \end{equation*}
    By looking at the diagram it is also
    obvious that the canonical injections
    in the coproduct commute with the
    semantic assignment, and since
    they send generating objects to
    generating objects, they are
    morphisms of $\PetriS$. They are clearly
    also transition-preserving.

    On morphisms, the monoidal structure
    acts like the coproduct of functors
    in the category of symmetric monoidal categories.
    Again, commutativity with
    assignment of semantics is guaranteed by
    coproduct laws.

    The monoidal unit is taken to be the
    category with no generating objects
    and no generating morphisms, having as
    semantics assignment the
    functor sending symmetries to symmetries
    and the monoidal unit to the
    monoidal unit.

    Associators and unitors are taken to
    be the ones of the coproduct.
\end{proof}
\begingroup
\def\thetheoremUnified{\ref{lem: characterizing coequalizers}}
\begin{lemma}
  Let $\CategoryC, \CategoryD$ be free strict
  symmetric monoidal categories, and let
  $F,G: \CategoryC \to \CategoryD$ be a couple
  of transition-preserving functors
  sending generating objects to generating objects.
  Denoting with $\GObj{\CategoryC}$
  and $\GMor{\CategoryC}$ the generating
  objects and morphisms, respectively, of $\CategoryC$, then the
  coequalizer $\CategoryE$ of $F,G$ is the
  following free symmetric strict monoidal category:
  \begin{itemize}
    \item Generating objects of $\CategoryE$ are
    $\GObj{\CategoryD}/\simeq_o$, where $\simeq_o$
    is the equivalence relation generated by
    \begin{equation*}
      \forall x \in \GObj{\CategoryC}.(Fx = Gx)
    \end{equation*}
    \item Given a generating morphism $f$ of $\CategoryC$, and
    denoting with $f_F$ and $f_G$ the generating morphisms
    of $\CategoryD$ such that $Ff = \sigma \Cp f_F \Cp \sigma'$ and
    $Gf = \varsigma \Cp f_G \Cp \varsigma'$, generating
    morphisms of $\CategoryE$ are
    $\GMor{\CategoryD}/\simeq_m$, where $\simeq_m$
    is the equivalence relation generated by
    \begin{equation*}
      \forall f \in \GMor{\CategoryC}.(f_F = f_G)
    \end{equation*}
  \end{itemize}
\end{lemma}
\addtocounter{theoremUnified}{-1}
\endgroup
\begin{proof}
    First, we prove that $\CategoryE$ is well-defined.
    On objects, we identify couples of
    objects of $\CategoryD$ which are
    hit by $F$ and $G$, respectively, for some
    $x$, and then take the transitive closure of
    this identification. The coequalizer
    $\coeq: \CategoryD \to \CategoryE$, on objects,
    is defined as sending each object
    to its equivalence class.

    On morphsims, we notice that since $F$ and $G$ are
    transition-preserving, then for each morphism
    generator $f$ in $\CategoryC$ it is
    $Ff = \sigma \Cp f_F \Cp \sigma'$ and
    $Gf = \varsigma \Cp f_G \Cp \varsigma'$, with the $\sigma$s
    and $\varsigma$s
    being symmetries and $f_F, f_G$ being morphism
    generators.
    This ensures that the equivalence relation
    generated by
    $\forall f \in \GMor{\CategoryC}.(f_F = f_G)$
    boils down to an equation between generating morphisms,
    and is thus well-defined. We use elements in
    $\GMor{\CategoryD}/\simeq_m$ to generate morphisms of
    $\CategoryE$, in particular by setting their sources and
    targets to be:
    \begin{equation*}
      \Source{[f]_{\simeq_m}} = \mathfrak{O}[\Source{f}]_{\simeq_o}
      \qquad
      \Target{[f]_{\simeq_m}} = \mathfrak{O}[\Target{f}]_{\simeq_o}
    \end{equation*}
    Where $\mathfrak{O}$ denotes any ordering function on
    the set of generating objects of $\CategoryE$.
    Of course, there are many different ways to
    choose $\mathfrak{O}$, which guarantee isomorphic results.
    This is not a problem since coequalizers are unique
    up to isomorphism. On the other hand, the well-ordering
    theorem guarantees that a choice of $\mathfrak{O}$ is
    always possible, and so we can always define
    sources and targets this way.

    Independence from the representative
    $f$ in the definition is guaranteed
    by the functoriality condition on objects:
    If $g \in [f]$ then by definition
    there is an $h$ in $\GMor{\CategoryC}$ such
    that $f = h_F$, $g=h_G$ and $Fh = \sigma \Cp h_F \Cp \sigma'$,
    $Gh = \varsigma \Cp h_G \Cp \varsigma'$.
    Now, functoriality of $F$ and $G$ implies
    $\Source{Fh} = F{\Source{h}}$ and
    $\Source{Gh} = G{\Source{h}}$, respectively.
    Applying monoidality and the fact that
    $F, G$ send generating objects to generating objects,
    the condition on the generating objects of $\CategoryE$
    implies that $\Source{Fh}$ and $\Source{Gh}$ are made
    of the same generating objects and only differ by a permutation.
    Considering that $Fh = \sigma \Cp h_F \Cp \sigma'$ and
    $Gh = \varsigma \Cp h_G \Cp \varsigma'$, also $\Source{h_F}$
    and $\Source{h_G}$ differ only by a permutation. Hence
    applying $\mathfrak{O}$ to $\Source{h_F}$ and $\Source{h_G}$
    leads to the same result. The
    same argument holds for the target of $f$, $g$.

    The functor $\coeq:\CategoryD \to \CategoryE$ is defined,
    on morphisms, by sending each generating morphism of
    $\CategoryD$ to its equivalence class, pre- and post-
    composed with the obvious symmetries to satisfy the
    functoriality condition on objects. By definition, $\coeq$
    is transition-preserving.

    If some other
    functor $H:\CategoryD \to \CategoryB$
    coequalizes $F$ and $G$, then it has at
    least to identify $Fx$ and $Gx$ for each $x$,
    and $Ff$ and $Gf$ for each generating
    morphism $f$.

    This guarantees that $\coeq^{-1} \Cp H$
    is a well defined function, both on objects
    and morphisms. Hence, we can
    send each object $A$ of $\CategoryE$
    to $(\coeq^{-1} \Cp H)A$ in $\CategoryB$, and
    each generating morphsim $[f]$ of $\CategoryE$
    to $(\coeq^{-1} \Cp H)[f]$ in $\CategoryB$,
    pre- and post- composing with the needed
    symmetries.
    This mapping satisfies the
    coequalizer commuting condition, while uniqueness
    follows from the fact that $\simeq_o$ and $\simeq_m$
    are the smallest equivalence relations identifying
    generating objects and morphisms as prescribed by
    $F$ and $G$.
\end{proof}
\begingroup
\def\thetheoremUnified{\ref{lem: identifying n places}}
\begin{lemma}
  Denote with $O_n$ the Petri net consisting of $n$ places
  and no transitions.
  Let $\NetSem{N}$ be an identification of $\NetSem{M}$ via
  $F: \NetSem{M} \to \NetSem{N}$ with witness $O_n,l,r$.

  Then there exist transition-preserving functors
  $l_1, r_1, \dots, l_n, r_n$ such that the following
  diagram commutes, where $\coeq_i$ is the coequalizer
  of $l_i, r_i$, and dashed arrows are obtained
  from the universal property of coequalizers:
  \begin{equation*}
    % \scalebox{0.8}{
    \begin{tikzpicture}[node distance=2cm,>=stealth',bend angle=45,auto,xscale=2.5cm]
        \node (1) at (0,0) {$\Free{M}$};
        \node (2) [right of=1] {$\Free{N}_1$};
        \node (3) [right of=2] {$\dots$};
        \node (n-1) [right= 1.5cm of 3] {$\Free{N}_{n-1}$};
        \node (n) [right= 1.5cm of n-1] {$\Free{N}$};
        \node (w) [left of=1] {$\Free{O_n}$};

        \node(1a) [above of=1] {$\Free{O_1}$};
        \node(2a) [above of=2] {$\Free{O_1}$};
        \node(na) [above of=n-1] {$\Free{O_1}$};

        \node(1b) [below of=1] {$\Semantics$};

        \draw[transform canvas={yshift=0.5ex},->] (w) to node[font=\tiny] {$l$} (1);
        \draw[transform canvas={yshift=-0.5ex},->](w) to node[font=\tiny,swap] {$r$} (1);
        \draw[transform canvas={xshift=0.5ex},->] (1a) to node[font=\tiny] {$l_1$} (1);
        \draw[transform canvas={xshift=-0.5ex},->](1a) to node[font=\tiny,swap] {$r_1$} (1);
        \draw[transform canvas={xshift=0.5ex},->] (2a) to node[font=\tiny] {$l_2$} (2);
        \draw[transform canvas={xshift=-0.5ex},->](2a) to node[font=\tiny,swap] {$r_2$} (2);
        \draw[transform canvas={xshift=0.5ex},->] (na) to node[font=\tiny] {$l_n$} (n-1);
        \draw[transform canvas={xshift=-0.5ex},->](na) to node[font=\tiny,swap] {$r_n$} (n-1);

        \draw[dotted, ->] (1) to node[font=\tiny] {$\coeq_1$} (2);
        \draw[dotted, ->] (2) to node[font=\tiny] {$\coeq_2$} (3);
        \draw[dotted, ->] (3) to node[font=\tiny] {$\coeq_{n-1}$} (n-1);
        \draw[dotted, ->] (n-1) to node[font=\tiny] {$\coeq_{n}$} (n);

        \draw[->] (1) to node[font=\tiny,swap] {$\FunM$} (1b);
        \draw[dashed, ->] (2) to node[font=\tiny] {} (1b);
        \draw[dashed, ->] (n-1) to node[font=\tiny] {} (1b);

        \draw[->] (n) to node[font=\tiny] {$\FunN$} (1b);
        \end{tikzpicture}
% }
  \end{equation*}
  Moreover, it is $F = \coeq_1 \Cp \dots \Cp \coeq_{n}$. In other
  words: The identification of a finite,
  arbitrary number of places can be
  performed in steps, where at each step
  no more than two places are identified.
\end{lemma}
\addtocounter{theoremUnified}{-1}
\endgroup
\begin{proof}
    We can prove the statement by induction on $n$.
    The base case, with $n=1$, is trivial, since
    the following diagram clearly commutes:
    \begin{equation*}
        \scalebox{0.8}{
    \begin{tikzpicture}[node distance=2cm,>=stealth',bend angle=45,auto]
        \node (1) at (0,0) {$\Free{M}$};
        \node (2) [right of=1] {$\Free{N}$};
        \node (w) [left of=1] {$\Free{O_1}$};

        \node(1a) [above of=1] {$\Free{O_1}$};

        \node(1b) [below of=1] {$\Semantics$};

        \draw[transform canvas={yshift=0.5ex},->] (w) to node {$l$} (1);
        \draw[transform canvas={yshift=-0.5ex},->](w) to node[swap] {$r$} (1);
        \draw[transform canvas={xshift=0.5ex},->] (1a) to node {$l$} (1);
        \draw[transform canvas={xshift=-0.5ex},->](1a) to node[swap] {$r$} (1);

        \draw[dotted, ->] (1) to node {$F$} (2);

        \draw[->] (1) to node[swap] {$\FunM$} (1b);
        \draw[->] (2) to node {$\FunN$} (1b);

        \end{tikzpicture}
}
    \end{equation*}
    For the induction step, first notice how
    $O_{n+1} = O_n + O_1$ (here the symbol $+$
    denotes the coproduct of nets) and $\Free{O_{n+1}} =
    \Free{O_n} + \Free{O_1}$ (here the symbol $+$
    denotes the coproduct of categories).
    Now consider
    the following commutative diagram:
    \begin{equation*}
        % \scalebox{0.8}{
    \begin{tikzpicture}[node distance=2cm,>=stealth',bend angle=45,auto]
        \node (1) at (0,0) {$\Free{M}$};
        \node (2) [right = 1.5cm of 1] {$\Free{N}_n$};
        \node (3) [right = 1.5cm of 2] {$\Free{N}_1$};
        \node (w) [left = 1.5cm of 1] {$\Free{O_{n+1}}$};
        \node (w') [left = 1.5cm of w] {$\Free{O_{n}}$};
        \node (w'') [above = 1.5cm of w'] {$\Free{O_{n+1}}$};
        \node (o) [above = 1.5cm of w] {$\Free{O_{1}}$};

        \node(2a) [above = 1.5cm of 2] {$\Free{O_1}$};
        \node(1b) [below = 1.5cm of 1] {$\Free{N}$};

        \draw[transform canvas={yshift=0.5ex},->] (w) to node {$l$} (1);
        \draw[transform canvas={yshift=-0.5ex},->](w) to node[swap] {$r$} (1);
        \draw[transform canvas={xshift=0.5ex},->] (2a) to node {$\iota_2 \Cp r \Cp \coeq_n$} (2);
        \draw[transform canvas={xshift=-0.5ex},->](2a) to node[swap] {$\iota_2 \Cp l \Cp \coeq_n$} (2);
        \draw[=, double, double distance=1mm] (o) to node {} (2a);
        % \draw[transform canvas={yshift=-0.5ex},-](o) to node[swap] {} (2a);
        \draw[dotted, ->] (1) to node {$\coeq_n$} (2);
        \draw[dotted, ->] (2) to node {$\coeq_1$} (3);
        \draw[dotted,->] (1) to node[swap] {$\coeq_{n+1}$} (1b);
        \draw[dashed, <->] (3) to node {$\exists !$} (1b);
        \draw[right hook-latex] (o) to node {$\iota_2$} (w);
        \draw[right hook-latex] (w') to node[swap] {$\iota_1$} (w);
        \draw[left hook-latex] (o) to node[swap] {$\iota_2$} (w'');
        \draw[right hook-latex] (w') to node {$\iota_1$} (w'');
        \draw[transform canvas={xshift=0.25ex, yshift=0.25ex},-] (w) to node {} (w'');
        \draw[transform canvas={xshift=-0.25ex, yshift=-0.25ex},-](w) to node[swap] {} (w'');
        \end{tikzpicture}
% }
    \end{equation*}
    Here, $\coeq_{n+1}$ is the coequalizer of $l, r$, while
    $\coeq_{n}$ is the coequalizer of $\iota_1 \Cp l$ and $\iota_1 \Cp r$.
    On the other hand, $\coeq_1$ is the coequalizer
    of $\iota_2 \Cp l \Cp \coeq_n$ and $\iota_2 \Cp r \Cp \coeq_n$.
    We want to prove the lemma for $\coeq_{n+1}$ assuming
    that it holds for $\coeq_n$.
    Looking at the diagram, we can infer that:
    \begin{itemize}
        \item $\coeq_n \Cp \coeq_1$ coequalizes the couple $\iota_1 \Cp l$
        and $\iota_1 \Cp r$ as well as the couple $\iota_2 \Cp l$
        and $\iota_2 \Cp r$:
        \begin{gather*}
            \iota_1 \Cp l \Cp \coeq_n \Cp \coeq_1 = (\iota_1 \Cp l) \Cp \coeq_n \Cp \coeq_1 =
            (\iota_1 \Cp r) \Cp \coeq_n \Cp \coeq_1 = \iota_1 \Cp r \Cp \coeq_n \Cp \coeq_1\\
            \iota_2 \Cp l \Cp \coeq_n \Cp \coeq_1 = (\iota_2 \Cp l \Cp \coeq_n) \Cp \coeq_1 =
            (\iota_2 \Cp r \Cp \coeq_n) \Cp \coeq_1 = \iota_2 \Cp r \Cp \coeq_n \Cp \coeq_1
        \end{gather*}
        Applying the universal property of
        coproducts, this means that
        \begin{gather*}
            l \Cp \coeq_n \Cp \coeq_1 = [\iota_1,\iota_2] \Cp l \Cp \coeq_n \Cp \coeq_1 =
            [\iota_1,\iota_2] \Cp r \Cp \coeq_n \Cp \coeq_1 = r \Cp \coeq_n \Cp \coeq_1
        \end{gather*}
        And hence, because of the universal
        property of coequalizers, there is
        a unique arrow from $\Free{N}_1$ to
        $\Free{N}$ making the diagram above commute;

        \item $\coeq_{n+1}$ coequalizes $l$ and $r$, and thus also
        $\iota_2 \Cp l \Cp \coeq_n$ and $\iota_2 \Cp l \Cp \coeq_n$. Again,
        applying the universal property of coequalizers
        we have that there exists a unique
        arrow from $\Free{N}$ to $\Free{N}_1$;

        \item The universal property of coequalizers forces
        these arrows to be one the inverse of the other. This
        means that $\Free{N}_1$ and $\Free{N}$ are isomorphic,
        and since coequalizers are unique up to isomorphism, that
        $\coeq_n \Cp \coeq_1$ is a coequalizer for $l,r$;

        \item Since $\coeq_n$ has $O_n$ as witness, we can
        apply the induction hypothesis, and infer the existence
        of $l_1, r_1, \dots, l_n, r_n$ and coequalizers $c_1,
        \dots, c_n$ such that $\coeq_n = c_1 \Cp \dots \Cp c_n$.

        \item Putting everything together, we have
        transition-preserving functors and coequalizers,
        respectively,
        \begin{equation*}
            l_1, r_1 \dots, l_n, r_n,(\iota_2 \Cp l \Cp \coeq_n), (\iota_2 \Cp l \Cp \coeq_n)
            \qquad
            c_1, \dots, c_n, \coeq_1
        \end{equation*}
        such that
        $c_1 \Cp \dots \Cp c_n \Cp \coeq_1 = \coeq_n \Cp \coeq_1 = \coeq_{n+1}$.
    \end{itemize}
    Commutativity with assignment of
    semantics is trivially satisfied applying the
    universal property of coequalizers.
    This concludes the proof.
\end{proof}
\begingroup
\def\thetheoremUnified{\ref{lem: uniqueness of erasing of generators}}
\begin{lemma}
  Given a net $\NetSem{N}$ and a subnet $N_w$,
  erasings of generators of $\NetSem{N}$ via
  $N_w$ are unique up to isomorphism.
\end{lemma}
\addtocounter{theoremUnified}{-1}
\endgroup
\begin{proof}
  We defer the proof to this lemma to the proof
  of Lemma~\ref{lem: characterization of erasing morphisms},
  of which this is a corollary.
\end{proof}
\begingroup
\def\thetheoremUnified{\ref{lem: decompositions are invariant for belonging}}
\begin{lemma}
    Let $\CategoryC$ be a FSSMC. Every morphism admits
    exactly one decomposition.
\end{lemma}
\addtocounter{theoremUnified}{-1}
\endgroup
\begin{proof}
    Since $\CategoryC$ is free, decompositions of $f$
    can differ only modulo symmetries and identities.
    This means that in each decomposition the
    same morphism generators are used, concluding the proof.
\end{proof}
\begingroup
\def\thetheoremUnified{\ref{lem: characterization of addition of generating morphisms}}
\begin{lemma}
  Let $\NetSem{M}$ be an addition of generating moprhisms to
  $\NetSem{K}$ via  $W, w$. If $W$ has finite places
  and transitions, then $\Free{M}$ is isomorphic to the
  FSSMC $\CategoryC$ generated as follows:
  \begin{itemize}
    \item $\CategoryC$ and $\Free{K}$ coincide on
    object generators (and hence on objects). Concisely:
    \begin{equation*}
      \GObj{\CategoryC} = \GObj{\Free{K}}
    \end{equation*}
    \item  The generating morphisms of $\CategoryC$
    are the disjoint union of the generating morphism
    of $\Free{K}$ and $\Free{W}$, where the domain
    and codomain of generators in $W$ are substituted
    with their images through $w$. Concisely:
    \begin{equation*}
      \GMor{\CategoryC} = \GMor{\Free{K}} \sqcup
      \Suchthat{\alpha: [\Source{\alpha}/\Source{w\alpha}] \to  [\Target{\alpha}/\Target{w\alpha}]}{\alpha \in \GMor{\Free{W}}}
    \end{equation*}
  \end{itemize}
  While $\FunM$ is the functor
  $\Free{M} \simeq \CategoryC \xrightarrow{F} \Semantics$,
  with $F$ defined as $\FunK$ on objects and on generating
  morphisms coming from $\Free{K}$, and as $w \Cp \FunK$
  on generating morphisms coming from $\Free{W}$.
\end{lemma}
\addtocounter{theoremUnified}{-1}
\endgroup
\begin{proof}
  It is clear that the coproduct of $\Free{K}$
  and $\Free{W}$ has the disjoint union
  of their object/morphism generators as
  object/morphism generators, respectively.

  Characterizing the pushout as the coequalizer
  of the coproduct as we already did in Section~\ref{sec: synchronizations},
  $\CategoryC$ is isomorphic to the coequalizer of
  $\Free{K} + \Free{W}$ through $\overline{W}$.
  Since $W$ is finite, the thesis follows by recursively applying
  Lemmas~\ref{lem: identifying two places} and~\ref{lem: identifying n places}.
\end{proof}
\begingroup
\def\thetheoremUnified{\ref{lem: characterization of erasing morphisms}}
\begin{lemma}
  Let $\NetSem{K}$ be an erasing of generators of
  $\NetSem{N}$ via $N_w$, with $N$ having
  finite places and transitions. Then $\Free{K}$
  is isomorphic to the FSSMC $\CategoryD$ generated
  as follows:
  \begin{itemize}
    \item $\CategoryD$ and $\Free{N}$ coincide on
    object generators (and hence on objects). Concisely:
    \begin{equation*}
      \GObj{\CategoryD} = \GObj{\Free{N}}
    \end{equation*}
    \item  The generating morphisms of $\CategoryD$
    are the generating morphisms of $\Free{N}$
    minus the generating morphisms of $\Free{N_w}$.
    Concisely:
    \begin{equation*}
      \GMor{\CategoryD} = \GMor{\Free{N}} /
      \Suchthat{\alpha}{\alpha \in \GMor{\Free{N_w}}}
    \end{equation*}
    Where $/$ denotes the set theoretic difference.
  \end{itemize}
  $\FunK$ is the functor
  $\Free{K} \simeq \CategoryD \xrightarrow{F} \Semantics$,
  with $F$ defined as $\FunN$ on objects and on generating
  morphisms.
\end{lemma}
\addtocounter{theoremUnified}{-1}
\endgroup
\begin{proof}
  By definition, it is clear that $\CategoryD$
  is an erasing of generators of $\NetSem{N}$
  via $N_w$. Now, let $\NetSem{K}$ be
  an erasing of generators of $\NetSem{N}$
  via $N_w$.
  Since $\overline{N_w}$ does not have any
  generating morphism, the generating morphisms of $\Free{N}$
  are just the disjoint union of the generating
  morphisms of $\Free{N_w}$ and the generating morphisms
  of $\Free{K}$. Then we have
  \begin{equation*}
   \GMor{\CategoryD} \sqcup \GMor{\Free{N_w}} \simeq
   \GMor{\Free{N}} \simeq
   \GMor{\Free{K}} \sqcup \GMor{\Free{N_w}}
  \end{equation*}
  From which it follows that the generating
  morphisms of $\Free{K}$ must be in bijection
  with the generating morphisms of $\CategoryD$
  since all the sets involved are finite.

  For generating objects, notice that the functor
  $\Free{\overline{N_w}} \to \Free{K}$ in the
  pushout square
  \begin{equation*}
    \scalebox{1}{
  \begin{tikzpicture}[node distance=2cm,>=stealth',bend angle=45,auto]
    \node (1) at (0,0) {$\Free{N}$};
    \node (2) at (4, 0) {$\Free{K}$};
    \node (3) at (0,2) {$\Free{N_w}$};
    \node (4) at (4,2) {$\Free{\overline{N_w}}$};

    \draw[left hook-latex] (4) to node[midway, above] {$\inj_{N_w}$} (3);
    \draw[->] (4) to (2);
    \draw[left hook-latex] (3) to node[midway, left] {$\sub_{N_w}$} (1);
    \draw[left hook-latex] (2) to node[midway, below] {$\sub_{K}$} (1);
  \end{tikzpicture}
}
  \end{equation*}
  must be injective on generating objects.
  To see this, denote with
  \begin{align*}
    \overline{\inj_{N_w}}:& \GObj{\Free{\overline{N_w}}} \to \GObj{\Free{N_w}}\\
    \overline{\sub_{N_w}}:& \GObj{\Free{N_w}} \to \GObj{\Free{N}}\\
    \overline{\sub_{K}}:&  \GObj{\Free{K}} \to \GObj{\Free{N}}
  \end{align*}
  The functions between sets obtained by suitably
  restricting $\inj_{N_w}$, $\sub_{N_w}$ and $\sub_K$,
  respectively. Clearly, all these functions are inclusions,
  so in particular injective.

  Now  denote with  $f: \GObj{\Free{\overline{N_w}}}
  \to \GObj{\Free{K}}$
  the restriction to generating objects
  of the functor $\Free{\overline{N_w}} \to \Free{K}$
  in the pushout square. Our claim amounts
  to prove that $f$ is injective.
  Suppose there are two functions
  $g,h: X \to \GObj{\overline{N_w}}$
  for some set $X$, such that $g \Cp f = h \Cp f$.
  Then it is
  $g \Cp f \Cp \overline{\sub_K} = h \Cp f \Cp \overline{\sub_k}$,
  and because of the commutativity of the square,
  $g \Cp \overline{\inj_{N_w}} \Cp \overline{\sub_{N_w}} =
  f \Cp \overline{\inj_{N_w}} \Cp \overline{\sub_{N_w}}$.
  But $\overline{\inj_{N_w}} \Cp \overline{\sub_{N_w}}$
  is injective being a composition of injective functions.
  Since in $\Set$ the injective functions are exactly the
  monomorphisms, it follows that $g = h$, proving that
  $f$ is mono, and hence injective.
  From this, we can write
  $\GObj{\Free{K}} = X \sqcup Y$
  for some sets $X, Y$,
  with $X \simeq \GObj{\Free{\overline{N_w}}}$.

  Now, notice that the objects in $Y$ are
  not identified with the generating objects of
  $\Free{N_w}$ by the pushout square.
  This implies that $\GObj{\Free{K}}$ is in bijection
  with $\GObj{\Free{N}}$, since the only way to
  obtain $\Free{N}$ in the pushout square is
  that $Y$ is in bijection with the generating
  objects of $\Free{N}$ that are not in $\Free{N_w}$.

  The claim on $\FunK$ is trivial.
 \end{proof}
\begingroup
\def\thetheoremUnified{\ref{lem: characterization of synchronization}}
\begin{lemma}
  Let $\NetSem{M}$ be a synchronization of
  $\NetSem{N}$ via $W, w$. If $W$
  has finite places and transitions, then $\Free{M}$
  is isomorphic to the FSSMC $\CategoryE$ generated
  as follows:
  \begin{itemize}
    \item $\CategoryE$ and $\Free{N}$ coincide on
    object generators (and hence on objects). Concisely:
    \begin{equation*}
      \GObj{\CategoryE} = \GObj{\Free{N}}
    \end{equation*}
    \item  The generating morphisms of $\CategoryE$
    are the generating morphisms of $\Free{N}$, to which
    are added the morphisms of $\Free{W}$ as in
    Lemma~\ref{lem: characterization of addition of generating morphisms},
    and to which are stripped away the generating
    morphisms of $\Free{N}$ that are in the decomposition of the image
    of generating morphisms of $\Free{W}$ through $w$.
    Concisely:
    \begin{align*}
      \GMor{\CategoryE} &=
      \GMor{\CategoryC} \sqcup \\
      & \sqcup \Suchthat{\alpha: [\Source{\alpha}/\Source{w\alpha}] \to  [\Target{\alpha}/\Target{w\alpha}]}{\alpha \in \GMor{\Free{W}}} / \\
      & /  \Suchthat{g \in \GMor{\Free{N}}}{\exists f \in \GMor{\Free{W}}.(g \in wf)}
    \end{align*}
    Where $/$ denotes the set theoretic difference.
  \end{itemize}
  $\FunM$ is the functor
  $\Free{M} \simeq \CategoryD \xrightarrow{F} \Semantics$,
  with $F$ defined as $\FunN$ on objects and on generating
  morphisms coming from $\Free{N}$, and as $w \Cp \FunN$
  on generating morphisms coming from $\Free{W}$.
\end{lemma}
\addtocounter{theoremUnified}{-1}
\endgroup
\begin{proof}
  Consider the pushout defining a synchronization:
  \begin{equation*}
    
  \end{equation*}
  $\Free{K}$ can be characterized by applying
  Lemma~\ref{lem: characterization of erasing morphisms}
  to $\NetSem{N}$ via $N_w$, while $\Free{M}$
  is obtained by applying
  Lemma~\ref{lem: characterization of addition of generating morphisms}
  to $\NetSem{K}$ via $W,w$. The thesis follows
  immediately by noticing that the set
  \begin{equation*}
    \Suchthat{g \in \GMor{\Free{N}}}{\exists f \in \GMor{\Free{W}}.(g \in wf)}
  \end{equation*}
  Is, by definition, the set of generating morphisms of $\Free{N_w}$.
\end{proof}

\end{document}